\documentclass[reqno, 12pt]{amsart}
\usepackage[margin=1in]{geometry}
\pdfoutput=1
\usepackage[T2A]{fontenc}
\usepackage[english]{babel}

\usepackage{amssymb,amsthm,amsfonts,amsmath,amscd,bbm,eucal,faktor,xfrac,braket,relsize,nccmath,mathtools,comment,verbatim,wasysym,graphicx,import,float,mathrsfs,color,xparse,chngcntr,mwe,enumitem,tikz,graphicx,wrapfig,lipsum,dsfont,etoolbox,stackrel,bm,upgreek}
\usepackage{url}

\usepackage[
backref=0,
backend=biber,
style=alphabetic,
citestyle=alphabetic
]{biblatex}

\usepackage[figurename=Figure.]{caption}

\usetikzlibrary{patterns}
\usetikzlibrary{shapes}

\usepackage[
  linktocpage=true,
]{hyperref}

\mathtoolsset{showonlyrefs,showmanualtags}

\allowdisplaybreaks

\definecolor{linkblue}{rgb}{0,0.2,0.6}
\definecolor{pictureblue}{rgb}{0,0,1}
\definecolor{pictureblack}{rgb}{0,0,0}
\definecolor{picturered}{rgb}{1,0,0}
\definecolor{picturegreen}{rgb}{0,1,0}

\newcommand{\quotid}[2]{{\raisebox{.3em}{$#1\!$}\bigl/\raisebox{-.3em}{$\!\left(#2\right)$}}}

\newcommand{\RR}{\mathbb{R}_{\geq 0}}
\newcommand{\ZZ}{\mathbb{Z}_{\geq 0}}

\newcommand\Z{\mathbb Z}
\newcommand\R{\mathbb R}

\newcommand\Pas{\mathbb{P}}

\newcommand\T{\mathcal{T}}

\newcommand\Ga{\Gamma}

\newcommand\ka{\varkappa}

\newcommand\la{\lambda}

\newcommand\wt{\widetilde}
\newcommand\wh{\widehat}

\newcommand\prodd{\prod\limits}
\newcommand\summ{\sum\limits}

\newcommand{\ol}{\overline}
\makeatletter
\newcommand*{\isomorphism}{%
  \mathrel{%
    \mathpalette\@isomorphism{}%
  }%
}
\newcommand*{\@isomorphism}[2]{%
  \sbox0{$#1\simeq$}%
  \sbox2{$#1\sim$}%
  \dimen@=\ht0 %
  \advance\dimen@ by -\ht2 %
  \sbox0{%
    \lower2.5\dimen@\hbox{%
      $\m@th#1\relbar\isomorphism@joinrel\rightarrow$%
    }%
  }%
  \rlap{%
    \hbox to \wd0{%
      \hfill\raise\dimen@\hbox{$\m@th#1\sim$}\hfill
    }%
  }%
  \copy0 %
}
\newcommand*{\isomorphism@joinrel}{%
  \mathrel{%
    \mkern-3.4mu %
    \mkern-1mu %
    \nonscript\mkern1mu %
  }%
}
\makeatother
\newcommand\iso{\raisebox{0.5ex}{$\, \isomorphism\, $}}
\newcommand\const{\operatorname{const}}

\newcommand\supp{\operatorname{supp}}
\newcommand{\spn}[2]{\operatorname{span}_{#1}\left(#2\right)}

\newcommand{\Res}{\operatorname{Res}}
\newcommand\Hom{\operatorname{Hom}}

\newcommand\Ext{\operatorname{Ext}}
\newcommand{\id}{\mathrm{id}}

\newcommand{\sat}{\operatorname{sat}}
\renewcommand\H{\mathcal{H}}
\DeclareMathOperator{\FH}{\mathcal{FH}}
\DeclareMathOperator{\exH}{\H_{ex}}
\DeclareMathOperator{\exFH}{\FH_{ex}}
\DeclareMathOperator{\exHo}{\H_{ex}^{\circ}}
\DeclareMathOperator{\exFHo}{\FH_{ex}^{\circ}}
\DeclareMathOperator{\K}{K_0}
\DeclareMathOperator{\KK}{K_0^+}
\newcommand{\MC}{\operatorname{MC}}

\newcommand{\lr}[1]{\left( #1 \right)}
\newcommand{\biglr}[1]{\Bigl( #1 \Bigr)}

\newcommand{\Inf}{\infty}

\renewcommand{\binom}[2]{\begin{pmatrix}
#1\\ 
#2
\end{pmatrix}
}
\newenvironment{speqn}
{
\begin{equation}
    \begin{gathered}
    }
    {
    \end{gathered}
\end{equation}
}
\newenvironment{speqn*}
{
\begin{equation*}
    \begin{gathered}
    }
    {
    \end{gathered}
\end{equation*}
}

\newtheorem{theorem}{Theorem}
\newtheorem{proposition}[theorem]{Proposition}
\newtheorem{corollary}[theorem]{Corollary}
\newtheorem{lemma}[theorem]{Lemma}

\newtheorem*{theorem*}{Theorem}
\newtheorem*{proposition*}{Proposition}
\newtheorem*{corollary*}{Corollary}

\theoremstyle{definition}
\newtheorem{definition}[theorem]{Definition}
\newtheorem{definition-proposition}{Definition-Proposition}
\newtheorem{remark}[theorem]{Remark}
\newtheorem{observation}[theorem]{Observation}
\newtheorem{notation-definition}[theorem]{Notation-Definition}
\newtheorem{example}[theorem]{Example}

\newtheorem*{notation}{Notation}
\newtheorem*{remark*}{Remark}
\newtheorem*{observation*}{Observation}
\newtheorem*{example*}{Example}

\numberwithin{definition-proposition}{section}
\numberwithin{theorem}{section}

\title{Semifinite harmonic functions on branching graphs}

\author{Nikita Safonkin}

\address{Skolkovo Institute of Science and Technology, Moscow, Russia \& National Research University Higher School of Economics, Moscow, Russia.}

\email{safonkin.nik@gmail.com}

\begin{document}

\begin{abstract}
    We study semifinite harmonic functions on arbitrary branching graphs. We give a detailed exposition of an algebraic method which allows one to classify semifinite indecomposable harmonic functions on some multiplicative branching graphs. This method was proposed by A. Wassermann in terms of operator algebras, while we rephrase, clarify, and simplify the main arguments, working only with combinatorial objects. 
    This work was inspired by the theory of traceable factor representations of the infinite symmetric group $S(\infty)$. 
\end{abstract}

    \maketitle
    \setcounter{tocdepth}{2}
    \tableofcontents

    \section{Introduction}
Classical character theory of finite and compact groups may be generalized to other classes of groups and algebras in different ways. For groups and $C^*$-algebras \textit{not} of type I the character theory deals \textit{not} with irreducible representations but with normal factor representations, i.e. homomorphisms to von Neumann algebras with a finite or semifinite trace. For AF-algebras one can reformulate the character theory in a combinatorial-algebraic language, speaking about non-negative harmonic functions on infinite graded graphs of a special type. Equivalently, one can treat these harmonic functions as central measures on the space of monotone paths in the graph. This approach was developed in works of A. M. Vershik and S. V. Kerov in the late of 70's --- early 80's. Harmonic functions that take only finite values lead to probability measures on the path space and correspond to factor representations of finite types I$_n$ and II$_1$. The connection between harmonic functions and normal factor representations motivates us to study the so-called \textit{semifinite} harmonic functions, which correspond to normal factor representations of type I$_{\infty}$ and II$_{\infty}$. These functions must take on a value $+\infty$ and satisfy some natural condition, see Definition \ref{semifinite def} below. 

A. M. Vershik and S. V. Kerov have classified semifinite harmonic functions on the Young and Kingman graphs, see \cites{versh_ker_81, kerov_vershik1990}. They have solved this problem with help of the so-called ergodic method, which involves an evaluation of a non-trivial limit. This method can be applied to any branching graph, but the main difficulty, which is not always easy to overcome, is to compute that limit. There is another approach developed by A. J. Wassermann. In his dissertation \cite{wassermann1981}, he suggested  to use a bijection between the faithfull factor representations of a primitive $C^*$-algebra and that of any its closed two-sided ideal, see \cite[p. 143, Theorem 7]{wassermann1981}. Another ingredient of Wassermann's method requires that the $K_0$-group of the corresponding AF-algebra admits a compatible ring structure, see \cite[p.146, Theorem 8]{wassermann1981}. Therefore, Wassermann's method is applicable to some multiplicative graphs only, for which it may be extremely useful. A. Wassermann has applied his method to determine all indecomposable semifinite harmonic functions on the Young graph and, thereby, he has proved the theorem of Vershik and Kerov without the ergodic method or any other complicated analytical computations.

The present paper contains a detailed exposition of Wassermann's method in terms of algebraic combinatorics unlike the original work \cite[chapter III, Section 6]{wassermann1981}, where the language of operator algebras was used. The combinatorial set-up allows us to clarify and simplify the main arguments of \cite[chapter III, Section 6]{wassermann1981}. Furthermore, we work with a generalization of branching graphs, namely, we consider branching graphs with formal non-negative multiplicities on edges. Crucial statements of Wassermann's method may be found in \cites[chapter III, Section 6]{wassermann1981}[]{strat_voic1975}[]{Bratteli1972} and \cite{vershik_Kerov83,versh_ker_85,kerov_vershik1990}. We prove them in a completely combinatorial way. These statements together with the original argument of A.Wassermann constitute a powerful method for the determination of indecomposable semifinite harmonic functions on those multiplicative graphs for which the limit from the Vershik-Kerov ergodic method turns out to be too complicated for an evaluation. The Macdonald graph, which corresponds to the simplest Pieri rule for the Macdonald symmetric functions, is a good example of such a graph. Using Wassermann's method one can obtain an exhaustive list of semifinite indecomposable harmonic functions on it, see Remark \ref{rem mac}. 

\subsection{Organization of the paper}
In Section \ref{2} we introduce graded graphs and discuss their ideals and coideals. In Section \ref{3} we introduce semifinite harmonic functions and prove some general facts about them. Section \ref{4} deals with semifinite harmonic functions on multiplicative branching graphs only. Section \ref{5} contains a combinatorial analog of an observation due to R. P. Boyer. In Appendix \ref{addendum} we discuss finite harmonic functions on the product of branching graphs.

\subsection{Acknowledgements}
The author is deeply grateful to Grigori Olshanski for many useful comments and stimulating discussions. I would also like to thank Pavel Nikitin for comments after reading the first draft of this paper. Supported in part by the Simons Foundation. I was partially supported by the HSE University Basic Research Program.

\section{Ideals and coideals of graded graphs}\label{2}
\phantom{}
In this section we recall main notions on branching graphs, ideals and coideals. 

\begin{definition}\label{def1}
    By a \textit{graded graph} we mean a pair $\lr{\Ga,\ka}$, where $\Ga$ is a graded set $\Ga=\bigsqcup\limits_{n\geq 0}\Ga_n$, $\Ga_n$ are finite sets and $\ka$ is a function $\Ga\times \Ga\rightarrow \RR$, that satisfies the following constraints:
    \begin{enumerate}[label=\theenumi)]
    \item\label{cond1} if $\la\in\Ga_n$ and $\mu\in\Ga_m$, then $\ka(\la,\mu)=0$ for $m-n\neq 1$. 
    
    \item\label{cond1.5} for any vertex $\la\in\Ga_n$ there exists $\mu\in\Ga_{n+1}$ with $\ka(\la,\mu)\neq 0$.
\end{enumerate}
Edges of the graded graph $(\Ga,\ka)$ are, by definition, pairs of vertices $(\la,\mu)$ with $\ka(\la,\mu)>0$. Then we may treat $\ka(\la,\mu)$ as a formal multiplicity of the edge.
\end{definition}

If $\la\in\Ga_n$, then the number $n$ is uniquely defined. We denote it by $|\la|$. We write $\la\nearrow \mu,$ if $|\mu|-|\la|=1$ and $\ka(\la,\mu)\neq 0$. In this case we say that \textit{there is an edge from $\la$ to $\mu$ of multiplicity $\ka(\la,\mu)$}.

Condition \ref{cond1} from Definition \ref{def1} means that we allow edges only between adjacent levels and condition \ref{cond1.5} means that each vertex must be connected by an edge with some vertex from the higher level.

\textit{A path} in a graded graph $\Ga$ is a (finite or infinite) sequence of vertices  $\la_1,\la_2,\la_3,\ldots$ such that $\la_i\nearrow \la_{i+1}$ for every $i$. We will write $\nu>\mu$ if $|\nu|>|\mu|$ and there is a path that connects $\mu$ and $\nu$. We write $\nu\geq \mu$, if $\nu=\mu$ or $\nu>\mu$. Relation $\geq$ turns $\Ga$ into a poset.

Let $\mu,\nu\in\Ga$ and $|\nu|-|\mu|=n\geq 1$. Then the following expression
\begin{speqn}\label{shiftdim}
\dim(\mu,\nu)=\summ_{\substack{\la_0,\ldots,\la_n\in\Ga:\\ \mu=\la_0\nearrow\la_1\nearrow\ldots\nearrow\la_{n-1}\nearrow\la_n=\nu}}\ka(\la_0,\la_1)\ka(\la_1,\la_2)\ldots\ka(\la_{n-1},\la_n).
\end{speqn}

is the "weighted" number of paths from $\mu$ to $\nu$. By definition we also set $\dim(\mu,\mu)=1$ and $\dim(\mu,\nu)=0$, if $\nu\not\geq\mu$. The function $\dim(\cdot,\cdot)\colon \Ga\times\Ga\rightarrow \R_{\geq 0}$ is called \textit{the shifted dimension}. Note that $\dim(\mu,\nu)=\ka(\mu,\nu)$, if $\mu\nearrow\nu$, and for $\mu\in\Ga_m$, $\nu\in\Ga_n$ and any $k$ such that $m\leq k\leq n$ we have

\begin{speqn}\label{dimension}
\dim(\mu,\nu)=\summ_{\la:\la\in\Ga_k}\dim(\mu,\la)\dim(\la,\nu).
\end{speqn}

\begin{definition}
    A \textit{branching graph} is defined as a graded graph $(\Ga,\ka)$ that satisfies the following conditions
    \begin{itemize}
        \item $\Ga_0=\{\diameter\}$ is a singleton,
        
        \item for any $\la\in\Ga_n$ with $n\geq 1$ there exists $\mu\in\Ga_{n-1}$ such that $\mu\nearrow\la$.
    \end{itemize}
\end{definition}

For a branching graph $(\Ga,\ka)$ we denote the expression $\dim(\diameter,\la)$ by $\dim(\la)$ and call it \textit{the dimension of} $\la$.

\begin{definition}
    A subset of vertices $I$ of a graded graph $\Ga$ is called an \textit{ideal}, if for any vertices $\la\in I$ and $\mu\in\Ga$ such that $\mu>\la$ we have $\mu\in I$. A subset $J\subset \Ga$ is called a \textit{coideal}, if for any vertices $\la\in J$ and $\mu\in\Ga$ such that $\mu<\la$ we have $\mu\in J$.
\end{definition}
	
\begin{remark}
Our terminology differs from the terminology of poset theory. Namely, our ideals and coideals are usually called \textit{filters} and \textit{ideals} respectively \cite{stanley1}.
\end{remark}
	
	There is a bijective correspondence $I\leftrightarrow\Ga\backslash I$ between ideals and coideals. Let $J$ be a coideal and $I=\Ga\backslash J$ be the corresponding ideal. Then the following conditions are equivalent:
	\begin{enumerate}[label=\theenumi)]
	    \item if $\left\{\mu \;\middle|\; \la\nearrow\mu\right\}\subset I$, then $\la\in I$
	    
	    \item for any $\la\in J$ there exists a vertex $\mu\in J$ such that $\la\nearrow \mu$.
	\end{enumerate}
	
	\begin{definition}
	An ideal $I$ and the corresponding coideal $J$ are called \textit{saturated}, if they satisfy the conditions above. A saturated ideal $I$ is called \textit{primitive}, if for any saturated ideals $I_1,I_2$ such that $I=I_1\cap I_2$ we have $I=I_1$ or $I=I_2$. A saturated coideal $J$ is called \textit{primitive}, if for any saturated coideals $J_1,J_2$ such that $J=J_1\cup J_2$ we have $J=J_1$ or $J=J_2$.
	\end{definition}
	
	The bijection $I\leftrightarrow\Ga\backslash I$ maps primitive saturated ideals to primitive saturated coideals and vice versa. We will also use the fact that ideals and saturated coideals are graded graphs themselves. 
	
	Let $\Ga$ be a branching graph. The space of infinite paths in $\Ga$ starting at $\diameter$ will be denoted by $\T(\Ga)$. To every path $\tau=\left(\diameter,\la_1\nearrow \la_2\nearrow \ldots\right)\in\T(\Ga)$ we associate the saturated primitive coideal $\Ga_{\tau}=\bigcup\limits_{n\geq 1}\{\la\in\Ga\mid \la\leq \la_n\}$.
	
	In the next proposition we give a combinatorial characterization of saturated primitive coideals of an arbitrary graded graph, see \cite{Bratteli1972}. Moreover, for branching graphs we describe all such coideals in terms of path coideals $\Ga_{\tau}$, see \cite{strat_voic1975} and \cite[p.129]{wassermann1981}.

	\begin{proposition}\label{co-ideal}
	\begin{enumerate}[leftmargin=5ex, label=\theenumi)]
	    \item A saturated coideal $J$ of a graded graph is primitive if and only if  for any two vertices $\la_1,\la_2\in J$ we can find a vertex $\mu\in J$ such that $\mu\geq\la_1,\la_2$.
	    
	    \item Every saturated primitive coideal of a branching graph is of the form $J=\Ga_{\tau}$ for some path $\tau\in \T(\Ga)$.
	\end{enumerate}
	\end{proposition}
	\begin{proof}	
	Let $J\subset \Ga$ be a saturated coideal. Suppose that there exist vertices $\la_1, \la_2\in J$, that do not possess a common majorant. Let us prove that $J$ may be presented as a union of two distinct proper saturated coideals. We need to introduce some notation. For any $\la\in J$ the subset of vertices of $J$, that lie above $\la$, will be denoted by $J^{\la}$, i.e. $J^{\la}=\{\mu\in J\mid \mu\geq \la\}$. For any subset $A\subset J$ we define $\downarrow A$ as the subset of vertices of $J$, that lie below some vertex of $A$, i.e. $\downarrow A=\{\mu\in J\mid \mu\leq \la,\ \text{for some}\ \la\in A\}$. Finally, for any ideal $I$ of $J$ the symbol $\sat\lr{I}$ stands for the minimal saturated ideal that contains $I$. In other words, $\sat\lr{I}$ consists of all the vertices of $I$ and all vertices $\la\in J$ such that $\left\{\mu \;\middle|\; \la\nearrow\mu\right\}\subset I$. With this notation in mind we set $J_1= \downarrow(J^{\la_1}), J_2= J\backslash \sat(J^{\la_1})$. It is not difficult to see that $J_1$ and $J_2$ are saturated coideals and their union coincides with $J$. Obviously, $\la_1\in J_1$ and $\la_1\notin J_2$. Next, we use the fact that vertices $\la_1$ and $\la_2$ do not possess a common majorant to show that $\la_2\in J_2$ and $\la_2\notin J_1$. Thus, $J_1$ and $J_2$ are proper distinct coideals of $J$.
	
	Now suppose that for any vertices $\la_1,\la_2\in J$ there exists $\mu\in J$ with $\mu\geq\la_1,\la_2$. We will show that $J=\Ga_{\tau}$ for some path $\tau\in\T(\Ga)$. Let us denote by $x_1,x_2,\ldots$ all the vertices of $J$ enumerated in any (fixed) order. Since $J$ is primitive, it follows that we can construct a sequence of vertices $y_1\leq y_2\leq \ldots$ of $J$ with the following properties
    \begin{equation}
	    \begin{aligned}
        &y_1=x_1,\\
        \end{aligned}
        \ \ \ \ \ \ \ 
        \begin{aligned}
        &y_2\geq y_1,\\
        &y_2\geq x_2,\\
        &y_2\in J,\\
        \end{aligned}
        \ \ \ \ \ \ \ 
        \begin{aligned}
        &y_3\geq y_2,\\
        &y_3\geq x_3,\\
        &y_3\in J,\\
        \end{aligned}
        \ \ \ \ \ \ \ 
        \begin{aligned}
        \ldots \\
        \ldots \\
        \ldots \\
        \end{aligned}
        \ \ \ \ \ \ \ 
        \begin{aligned}
        &y_{n}\geq y_{n-1}\\
        &y_{n}\geq x_n\\
        &y_{n}\in J,
        \end{aligned}
        \ \ \ \ \ \ \ 
        \begin{aligned}
        \ldots \\
        \ldots \\
        \ldots \\
        \end{aligned}
    \end{equation}
Let $\tau\in\T(\Ga)$ be any path that goes through the vertices $y_1,y_2,\ldots$. Obviously,  $J=\Ga_{\tau}$.
\end{proof}

\begin{remark}
One can formulate an obvious analog of the second part of Proposition \ref{co-ideal} for arbitrary graded graphs, but this is of no particular importance to us.
\end{remark}
	
\begin{definition}
A graded graph $\Ga$ is called \textit{primitive} if it is primitive as a coideal, i.e. for any vertices $\la_1,\la_2\in \Ga$ there exists a vertex $\mu\in \Ga$ such that $\mu\geq\la_1,\la_2$.
\end{definition}

\section{Semifinite harmonic functions}\label{3}
\begin{definition}\label{harm def}
Let $(\Ga,\ka)$ be a graded graph. A function $\varphi\colon \Ga\rightarrow \R_{\geq 0}\cup\{+\Inf\}$ is called \emph{harmonic}, if it enjoys the following property:
\begin{speqn}
 \varphi(\la)=\summ_{\mu:\la\nearrow \mu}\ka(\la,\mu)\varphi(\mu),\ \ \forall \la\in\Ga.
\end{speqn}
\end{definition}

Throughout the paper we use the following conventions: 
\begin{itemize}
    \item $x+\lr{+\Inf}=+\Inf,\ \text{for any}\ x\in\R$,
    
    \item $\lr{+\Inf}+\lr{+\Inf}=+\Inf$,
    
    \item $0\cdot \lr{+\Inf}=0$.
\end{itemize} 

\begin{definition}
The set of all vertices $\la\in\Ga$ with $\varphi(\la)<+\Inf$ is called the \textit{finiteness ideal} of $\varphi$. We denote the \textit{zero ideal} $\left\{\la\in\Ga\; \middle|\; \varphi(\la)=0\right\}$ by $\ker{\varphi}$ and the \textit{support} $\{\la\in\Ga\mid \varphi(\la)>0\}$ by $\supp{\varphi}$.
\end{definition}

Note that the zero set $\ker(\varphi)$ is a saturated ideal and $\supp(\varphi)$ is a saturated coideal of $\Ga$ and $\ker(\varphi)\cup\supp(\varphi)=\Ga$. Furthermore, we can restrict $\varphi$ to any ideal or saturated coideal that contains $\supp(\varphi)$. The restriction is a harmonic function on that ideal or coideal respectively.

The symbol $\K(\Ga)$ stands for the $\R$-vector space spanned by the vertices of $\Ga$ subject to the following relations 
$$\la=\summ_{\mu:\la\nearrow \mu}\ka(\la,\mu)\cdot \mu,\ \ \ \forall \la\in\Ga.$$ 
The symbol $\KK(\Ga)$ denotes the positive cone in $\K(\Ga)$, generated by the vertices of $\Ga$, i.e. $\KK(\Ga)=\spn{\RR}{\la\mid\la\in\Ga}$. The partial order, defined by the cone $\KK(\Ga)$, is denoted by $\geq_K$. That is $a\geq_K b\iff a-b\in \KK(\Ga)$. For instance, if $\la\leq \mu$, then $\la\geq_K \dim(\la,\mu)\cdot\mu$. 

\begin{remark}
Notation $\K(\Ga)$ is motivated by the following fact. If all formal multiplicities of edges are integer numbers, then the vector space $\K(\Ga)$ can be identified with the Grothendieck $\K$-group of the corresponding AF-algebra. Under such a bijection the cone $\KK(\Ga)$ gets identified with the cone of true modules \cite[Theorem 13 on page 32]{versh_ker_85}. 
\end{remark}

\begin{observation}\label{rem cone}
If $b\in \KK(\Ga)$ and $b\leq_K \la$ then $b$ has the form $$b=\summ_{\mu\colon |\mu|=N}b_{\mu}\mu$$
for some $N$ and some real numbers $b_{\mu}$ subject to the following constraints $0 \leq b_{\mu}\leq \dim(\la,\mu)$. In particular, $b_{\mu}=0$, if $\mu\ngeq\la$.
\end{observation}

The $\RR$-linear map $\KK(\Ga)\rightarrow \RR\cup\{+\Inf\}$, defined by a harmonic function $\varphi$, will be denoted by the same symbol $\varphi$. Note that this map is monotone in the sense of the partial order. Namely, if $a\geq_K b$, then $\varphi(a)\geq \varphi(b)$.

\begin{definition}\label{semifinite def}
A harmonic function $\varphi$ is called \emph{semifinite}, if it is not finite and the map $\varphi\colon\KK(\Ga)\rightarrow \RR\cup\{+\Inf\}$ enjoys the following property
\begin{speqn}\label{lower semicont}
\varphi(a)=\sup_{\substack{b\in\KK(\Ga)\colon  b\leq_K a,\\ \varphi(b)<+\Inf}}\varphi(b),\ \ \ \ \ \forall a\in\KK(\Ga). 
\end{speqn}
\end{definition}

If $\varphi(a)<+\Inf$, then condition \eqref{lower semicont} becomes the trivial identity $\varphi(a)=\varphi(a)$.

Condition \eqref{lower semicont} arises in the theory of operator algebras in a natural way \cite[Definition 1.8]{boyer}.

\begin{remark}\label{remark 321}
 A harmonic function $\varphi$ is semifinite if and only if  there exists an element $a\in\KK(\Ga)$ with $\varphi(a)=+\Inf$ and for any such $a$ we can find a sequence $\{a_n\}_{n\geq 1}\subset \KK(\Ga)$ such that
 \begin{itemize}
     \item $a_n\leq_K a$,
     
     \item $\varphi(a_n)<+\Inf$,
     
     \item $\lim\limits_{n\to+\Inf}\varphi(a_n)=+\Inf$.
 \end{itemize} 
 We will call this $\{a_n\}_{n\geq 1}$ \textit{an approximating sequence}.
\end{remark}

\begin{proposition}\label{prop semifin redef}
A harmonic function $\varphi$ is semifinite if and only if  it is not finite and for any vertex $\la\in\Ga$ the following equality holds
\begin{speqn}\label{main equation}
   \varphi(\la)=\lim_{N\to\Inf}\summ_{\substack{\mu\colon\mu\geq \la,\ |\mu|=N\\ 0<\varphi(\mu)<+\Inf}}\dim(\la,\mu)\varphi(\mu).
\end{speqn}
\end{proposition}
\begin{proof}
If equality \eqref{main equation} is fulfilled, then $\varphi$ is semifinite, since prelimit sums give us an approximating sequence. If $\varphi$ is semifinite and $\varphi(\la)<+\infty$ then equality \eqref{main equation} is a trivial consequence of Definition \ref{harm def}. If $\varphi(\la)=+\infty$, then we can find an approximating sequence and Observation \ref{rem cone} yields that the prelimit expression is unbounded in $N$. We are left to prove that the limit exists. In fact, we show that the prelimit sequence is non-decreasing in $N$. Let us denote the prelimit expression by $\psi_N$. 

Next, the function 
$$\phi(\la)=\begin{cases}
   \varphi(\la),\ \text{if}\ 0<\varphi(\la)<+\Inf\\
   0\ \text{otherwise}
   \end{cases}$$ 
is \textit{subharmonic}: $$\phi(\la)\leq\summ_{\mu\colon\la\nearrow \mu}\ka\lr{\la,\mu}\phi(\mu).$$ 
Then from 
$$\psi_N=\summ_{\mu\colon |\mu|=N}\dim(\la,\mu)\phi(\mu)$$ 
and equality \eqref{dimension} it follows that $\psi_1\leq \psi_2\leq \psi_3\leq \ldots.$
\end{proof}
\begin{corollary}\label{useful corollary}
If $\varphi$ is a semifinite harmonic function on a graded graph $\Ga$, then for any vertex $\la\in\Ga$ with $\varphi(\la)=+\Inf$ there exists a vertex $\mu\geq \la$ such that $0<\varphi(\mu)<+\Inf$.
\end{corollary}

\begin{remark}\label{67}
Let $\{c_{\mu}\}_{\mu\in\Ga}$ be a tuple of non-negative real "numbers" $c_{\mu}\in\RR\cup\{+\Inf\}$ such that for every vertex $\la\in\Ga$ there exists the limit $\lim\limits_{N\to\Inf}\summ_{\mu\in\Ga_N}\dim(\la,\mu)c_{\mu}$, which may be infinite. For instance, we may take $c_{\mu}=\psi(\mu)$, where $\psi$ is a  \textit{subharmonic} function: $\psi(\la)\leq \summ_{\mu:\la\nearrow \mu}\ka(\la,\mu)\psi(\mu)$. Then the function 
$$\ol{c}(\la)=\!\lim\limits_{N\to\Inf}\!\summ_{\mu\in\Ga_N}\!\!\dim(\la,\mu)c_{\mu}$$ is harmonic, cf. \cite[p.4]{kerov_okounkov_olshanski}, see also \cite[formula (47)]{kerov_vershik1990}.
\end{remark}

\begin{definition}\label{def indec}
A semifinite harmonic function $\varphi$ is called \textit{indecomposable}, if for any finite or semifinite harmonic function $\varphi'$ which does not vanish identically on the finiteness ideal of $\varphi$ and satisfies the inequality $\varphi'\leq \varphi$ we have $\varphi'=\const\cdot\varphi$ on the finiteness ideal of $\varphi$.
\end{definition}

At the first glance the finiteness ideal of $\varphi'$ might be bigger than that of $\varphi$, but the next remark shows that this is not the case.                

\begin{remark}\label{rem indec}
If $\varphi$ and $\varphi'$ from Definition \ref{def indec} are proportional on the finiteness ideal of $\varphi$, then they are proportional on the whole graph $\Ga$. Indeed, by virtue of Proposition \ref{prop semifin redef} we may write
\begin{multline}
   \varphi(\la)=\const^{-1}\cdot\lim_{N\to\Inf}\summ_{\substack{\mu\colon\mu\geq \la,\ |\mu|=N\\ 0<\varphi(\mu)<+\Inf}}\dim(\la,\mu)\varphi'(\mu)\leq\\
   \leq \const^{-1}\cdot\lim_{N\to\Inf}\summ_{\substack{\mu\colon\mu\geq \la,\ |\mu|=N\\ 0<\varphi'(\mu)<+\Inf}}\dim(\la,\mu)\varphi'(\mu)=\const^{-1}\cdot\varphi'(\la).
\end{multline}

Thus, $\varphi'\leq\varphi\leq \const^{-1}\cdot\varphi'$ and finitiness ideals of $\varphi$ and $\varphi'$ coincide.
\end{remark}

\begin{notation}
The set of all indecomposable finite (not identically zero) and semifinite harmonic functions on a graded graph $\Ga$ is denoted by $\exH(\Ga)$. The subset of $\exH(\Ga)$ consisting of strictly positive functions is denoted by $\exHo(\Ga)$.
\end{notation}

\begin{lemma}\label{useful lemma}
Let $I$ be an ideal of a graded graph $\Ga$. Assume that $\varphi\in\exH(\Ga)$ does not vanish on $I$ identically. Then the following equality holds
\begin{speqn}\label{useful eq}
   \varphi(\la)=\lim\limits_{N\to\infty}\summ_{\substack{\mu:\mu\in I\\ |\mu|=N}}\dim(\la,\mu)\varphi(\mu),\ \ \la\in\Ga.
\end{speqn}
Moreover, for any element $a\in\KK(\Ga)$ we have $\varphi(a)=\sup\limits_{\substack{b\in\KK(I)\colon  b\leq_K a,\\ \varphi(b)<+\Inf}}\varphi(b)$.
\end{lemma}

\begin{remark}\label{rem after useful lemma}
If we omit the assumption that $\varphi$ is indecomposable, then the equality above should be replaced by the inequality $$\varphi(\la)\geq \lim\limits_{N\to\infty}\summ_{\substack{\mu:\mu\in I\\ |\mu|=N}}\dim(\la,\mu)\varphi(\mu).$$
\end{remark}
\begin{proof}[Proof of Lemma \ref{useful lemma}]
First of all, we remark that there exists a vertex $\nu\in I$ such that $0<\varphi(\nu)<+\infty$. Indeed, $\varphi$ does not equal zero identically on $I$, hence we can find a vertex $\nu'\in I$ such that $\varphi(\nu')>0$. If $\varphi(\nu')=+\infty$, then by Corollary \ref{useful corollary} we can find another vertex $\nu>\nu'$ with $0<\varphi(\nu)<+\infty$, which necessarily lies in $I$. 

Note that the function 
$$\phi(\la)=\begin{cases}
   \varphi(\la),\ \text{if}\ \la\in I,\\
   0\ \text{otherwise}
   \end{cases}$$
is subharmonic on $\Ga$. Then by Remark \ref{67} the right-hand side of \eqref{useful eq} defines a harmonic function on $\Ga$. From Observation \ref{rem cone} and Remark \ref{remark 321} it follows that the restriction of $\varphi$ to the ideal $I$ is a finite or semifinite harmonic function on $I$. Then the harmonic function on $\Ga$ defined by the right-hand side of \eqref{useful eq} is finite or semifinite as well. Next, by the very definition of harmonic functions, the prelimit expression is majorized by $\varphi$ for any $N$. Then the harmonic function that is defined as the limit $N\to +\infty$ is also majorized by $\varphi$. Finally, indecomposibility of $\varphi$ implies that $\varphi$ and the right-hand side of \eqref{useful eq} are proportional, but they coincide on the ideal $I$. Thus, they coincide on the whole graph $\Ga$, since there exists $\nu\in I$ with $0<\varphi(\nu)<+\infty$.
\end{proof}

Now we are ready to prove the most crucial statement of Wassermann's method. The following theorem is a combinatorial analog of a result, which is well known in the context of $C^*$-algebras, see \cites[Theorem 7 p.143, Corollary p.144 ]{wassermann1981} and \cite[II.6.1.6 p.102]{blackadar}. 
    
\begin{theorem}\label{theorem ext}
Let $I$ be an ideal of a graded graph $\Ga$.
\begin{enumerate}[leftmargin=3ex,label=\theenumi)]
\item There is a bijective correspondence between $\left\{\varphi\in\exH(\Ga)\colon \varphi\!\left.\right|_{I}\neq 0\right\}$ and $\exH(I)$, defined by the following mutually inverse maps
\begin{align}
    &\Res^{\Ga}_{I}:\left\{\varphi\in\exH(\Ga)\colon \varphi\!\left.\right|_{I}\neq 0\right\}\rightarrow \exH(I),\ \ \varphi\mapsto \varphi{\big|}_I, \\
    &\Ext^{\Ga}_{I}:\exH(I)\rightarrow \left\{\varphi\in\exH(\Ga)\colon \varphi\!\left.\right|_{I}\neq 0\right\},\ \ \varphi(\cdot)\mapsto \lim\limits_{N\to\infty}\summ_{\substack{\mu:\mu\in I\\ |\mu|=N}}\dim(\cdot,\mu)\varphi(\mu).
\end{align}
Furthermore, for any element $a\in\KK(\Ga)$ we have $\Ext^{\Ga}_{I}(\varphi)(a)=\sup\limits_{\substack{b\in\KK(I)\colon  b\leq_K a,\\ \varphi(b)<+\Inf}}\varphi(b).$

\item If $\Ga$ is a primitive graded graph, then the bijection above preserves strictly positive harmonic functions $\exHo\lr{I}\longleftrightarrow \exHo\lr{\Ga}$. 
\end{enumerate}
\end{theorem}
\begin{proof}
Suppose that $\varphi\in\exH(\Ga)$ and $\varphi\!\left.\right|_{I}\neq 0$. Then from Observation \ref{rem cone} and Remark \ref{remark 321} it follows that $\Res^{\Ga}_{I}\lr{\varphi}=\varphi\left.\right|_{I}$ is a finite or semifinite harmonic function on $I$. Lemma \ref{useful lemma} implies that $\Res^{\Ga}_{I}\lr{\varphi}$ is indecomposable.

Now let $\varphi\in\exH(I)$. From the proof of Proposition \ref{prop semifin redef} it follows that the limit from the definition of $\Ext^{\Ga}_{I}$ exists and $\Ext^{\Ga}_{I}\lr{\varphi}$ is a finite or semifinite harmonic function on $\Ga$. Note that $\Ext^{\Ga}_{I}\lr{\varphi}$ is strictly positive for $\varphi\in\exHo(I)$ because of the following simple fact, which holds for any primitive graded graph. For any vertex $\la\in\Ga$ there exists a vertex $\mu\in I$ such that $\mu\geq \la$.

Let us show that the harmonic function $\Ext^{\Ga}_{I}\lr{\varphi}$ is indecomposable for any $\varphi\in\exH(I)$. Suppose that $\Ext^{\Ga}_{I}\lr{\varphi}\geq \psi$, for some $\psi$, that does not vanish on the finiteness ideal of $\Ext^{\Ga}_{I}\lr{\varphi}$ identically. We denote that ideal by $\wt{I}$. The finiteness ideal of $\varphi$ is denoted by $I^{\varphi}$. Let us introduce more notation: $\psi_1=\psi\left.\right|_{\wt{I}}$ and $\psi_2=\Ext^{\wt{I}}_{I}(\varphi)-\psi_1$. Then $\psi_1$ and $\psi_2$ are finite harmonic functions on $\wt{I}$. Note that $\Ext^{\wt{I}}_{I}\lr{\varphi}=\Ext^{\wt{I}}_{I\cap I^{\varphi}}\lr{\varphi}$. On the one hand, we have $\Ext^{\wt{I}}_{I\cap I^{\varphi}}\lr{\varphi}=\psi_1+\psi_2$. On the other hand, $\varphi=\psi_1+\psi_2$ on $I\cap I^{\varphi}$, hence 
$$\Ext^{\wt{I}}_{I\cap I^{\varphi}}\lr{\varphi}=\Ext^{\wt{I}}_{I\cap I^{\varphi}}(\psi_1)+\Ext^{\wt{I}}_{I\cap I^{\varphi}}(\psi_2)\leq \psi_1+\psi_2,$$
where the last inequality follows from Remark \ref{rem after useful lemma}. Therefore, $\psi_1=\Ext^{\wt{I}}_{I\cap I^{\varphi}}(\psi_1)$ and $\psi_2=\Ext^{\wt{I}}_{I\cap I^{\varphi}}(\psi_2)$. Let us rewrite the first equality in the form $\psi\left.\right|_{\wt{I}}=\Ext^{\wt{I}}_{I\cap I^{\varphi}}(\psi)$. Then we see that the function $\psi\left.\right|_{I^{\varphi}}$ is not equal to zero identically. Now indecomposability of $\varphi$ yields that $\varphi$ and $\psi$ are proportional on $I^{\varphi}$. Thus, from $\psi\left.\right|_{\wt{I}}=\Ext^{\wt{I}}_{I\cap I^{\varphi}}\lr{\psi}$ and $\Ext^{\wt{I}}_{I}\lr{\varphi}=\Ext^{\wt{I}}_{I\cap I^{\varphi}}\lr{\varphi}$ it follows that $\Ext^{\Ga}_{I}\lr{\varphi}$ and $\psi$ are proportional on $\wt{I}$.

Therefore, maps $\Res^{\Ga}_{I}$ and $\Ext^{\Ga}_{I}$ are well defined and the following identity $\Res^{\Ga}_{I}\circ\Ext^{\Ga}_{I}=\id$ holds. The remaining identity $\Ext^{\Ga}_{I}\circ\Res^{\Ga}_{I}=\id$ immediately follows from Lemma \ref{useful lemma}. 
\end{proof}
\begin{remark}
Let $I_1\subset I_2$ be ideals of $\Ga$. Then $\Ext_{I_2}^{\Ga}\circ \Ext_{I_1}^{I_2}=\Ext_{I_1}^{\Ga}$.
\end{remark}

\begin{proposition}\label{prop prim coideal}\cite[p.35 Lemma 12]{versh_ker_85}
	Let $\Ga$ be a graded graph. If $\varphi\in\exH(\Ga)$, then the support $\supp(\varphi)$ is a primitive coideal.
\end{proposition}
\begin{proof}
	Let $\la_1,\la_2\in \supp(\varphi)$. Then Lemma \ref{useful lemma} yields $$\varphi(\la_2)=\lim\limits_{N\to\infty}\summ_{\substack{\mu:\mu\in \Ga^{\la_1}\\ |\mu|=N}}\dim(\la_2,\mu)\varphi(\mu),$$ 
	where $\Ga^{\la_1}=\{\nu\in\Ga\mid\nu\geq\la_1\}$. Then the inequality $\varphi(\la_2)>0$ implies that there exists a vertex $\mu$ such that $\mu\geq \la_1,\la_2$ and $\varphi(\mu)\neq 0$. Thus, by virtue of Proposition \ref{co-ideal} the coideal $\supp(\varphi)$ is primitive.
\end{proof}

\section{Multiplicative branching graphs}\label{4}

In this section we recall some basic notions related to multiplicative branching graphs \cite{kerov_book, versh_ker_85}. For such graphs we prove a theorem, which states that some multiplicative branching graphs admit no strictly positive semifinite indecomposable harmonic functions \cite[Theorem 8 p.146]{wassermann1981}. We call this theorem Wassermann's forbidding theorem. We also prove a semifinite analog of the Vershik-Kerov ring theorem \cite[Theorem p.144]{vershik_Kerov83}. 

\begin{definition}\label{multiplicative def}\cite[p.40]{versh_ker_85}
	A branching graph $\Ga$ is called \textit{multiplicative}, if there exists an associative $\Z_{\geq 0}$-graded $\R$-algebra $A=\bigoplus\limits_{n\geq0} A_n$, $A_0=\R$ with a distinguished basis of homogeneous elements $\{a_{\la}\}_{\la\in\Ga}$, that satisfy the following conditions
	\begin{enumerate}[label=\theenumi)]
	    \item $\deg{a_{\la}}=|\la|$
	    
	    \item $a_{\diameter}$ is the identity in $A$
	    
	    \item\label{multiplicative def3} For $\wh{a}=\summ_{\nu\in\Ga_1}\ka(\diameter,\nu)a_{\nu}$ and any vertex $\la\in\Ga$ we have $\wh{a}\cdot a_{\la}=\summ_{\mu:\la\nearrow\mu}\ka(\la,\mu)a_{\mu}$.
	\end{enumerate}
	
	Moreover, we assume that the structure constants of $A$ with respect to the basis $\{a_{\la}\}_{\la\in\Ga}$ are non-negative. 
\end{definition}

Let $\lr{\Ga,\ka}$ be the multiplicative graph that is related to an algebra $A$ and a basis $\{a_{\la}\}_{\la\in\Ga}$. We denote the quotient algebra $\quotid{A}{\wh{a}-1}$ by $R$, the canonical homomorphism $A\twoheadrightarrow R$ by $[\, \cdot \,]$ and the positive cone in $R$, consisting of all elements that can be written in the form $\summ_{\la\in\Ga_n}c_{\la}[a_{\la}]$ for a large enough $n$ and some $c_{\la}\geq 0$, by $R^+$. The correspondence $[\la]\mapsto [a_{\la}]$ defines an isomorphism of $\R$-vector spaces $\K(\Ga)\iso R $. The image of the cone $\KK(\Ga)\subset\K(\Ga)$ under this map coincides with $R^+$. 

Consider the positive cone $A^+\subset A$, consisting of all elements of $A$, that can be written as a linear combination of basis elements $a_{\la}$ with non-negative coefficients. For any semifinite harmonic function $\varphi\in\H(\Ga)$ we may speak about the $\R_{\geq0}$-linear map $\varphi\colon A^+\rightarrow \R_{\geq0}\cup\{+\Inf\}$. 

Let us now formulate the Vershik-Kerov ring theorem \cite[Theorem p.134]{vershik_Kerov83}, see also \cite[Proposition 8.4]{gnedin_olsh2006}.

\begin{definition}
A harmonic function $\varphi$ on a branching graph $\Ga$ is called \textit{normalized} if $\varphi(\diameter)=1$. 
\end{definition}

\begin{theorem}[Vershik-Kerov Ring Theorem]\cite[Theorem p.134]{vershik_Kerov83}\label{vershik kerov ring theorem}
A finite normalized harmonic function $\varphi$ on the multiplicative branching graph $\Ga$ is indecomposable if and only if the corresponding functional on $A$ is multiplicative: $\varphi\lr{a\cdot b}=\varphi\lr{a}\cdot \varphi\lr{b}\ \forall a,b\in A$.
\end{theorem}

The following semifinite analog of the ring theorem holds.

\begin{theorem}\cite[Theorem p.144]{vershik_Kerov83}\label{multiplicativity theorem}
For any semifinite indecomposable harmonic function $\varphi$ on the multiplicative branching graph $\Ga$ there exists a finite normalized indecomposable harmonic function $\psi$, such that $\varphi(a\cdot b)=\psi(a)\cdot\varphi(b)$ for any $a,b\in A^+$ with $\varphi(b)<+\Inf$.
\end{theorem}
\begin{proof}
Note that 
$$\lr{\wh{a}}^n=\summ_{\nu\colon\nu\in\Ga_n}\dim(\nu)\cdot a_{\nu}.$$ 
Then $\varphi\biglr{\lr{\wh{a}}^n a_{\mu}}=\varphi\lr{a_{\mu}}\geq \dim(\la)\varphi\lr{a_{\la}\cdot a_{\mu}}$ and $\varphi^{\la}\lr{\mu}=\varphi\lr{a_{\la}a_{\mu}}$ is a finite harmonic function on the finiteness ideal of $\varphi$. Since the restriction of $\varphi$ to its finiteness ideal is an indecomposable harmonic function (see Lemma \ref{useful lemma}) it follows that there exists $c_{\la}\in \R_{\geq 0}$ such that $\varphi\lr{a_{\mu}\cdot a_{\la}}=c_{\la}\varphi(a_{\mu})$. We set $\psi(\la)=c_{\la}$ by definition. One can check that $\psi$ is a harmonic function and that the functional on $A$ defined by $\psi$ is multiplicative. Then the Vershik-Kerov ring theorem implies that $\psi$ is indecomposable. 
\end{proof}

From Theorem \ref{multiplicativity theorem} it follows that the subspace $I=\spn{\R}{a_{\la}\mid\la\colon\varphi(\la)<+\Inf}\subset A$ is an ideal for any semifinite indecomposable harmonic function $\varphi$. However, the proof shows that this is true for an arbitrary harmonic function $\varphi$ without any additional assumptions.

The following theorem imposes some restrictions on multiplicative graphs that possess strictly positive indecomposable semifinite harmonic functions, \cite[Theorem 8 p.146]{wassermann1981}.
	
\begin{theorem}[Wassermann's forbidding theorem]\label{theorem mult graph}
If $a_{\la}a_{\mu}\neq0$ for any $\la,\mu\in\Ga$, then the graph $\Ga$ admits no strictly positive semifinite indecomposable harmonic functions.
\end{theorem}
\begin{proof}
Let $\varphi$ be a strictly positive indecomposable semifinite harmonic function. The argument given at the beginning of the proof of Theorem \ref{multiplicativity theorem} shows that $\varphi^{\mu}$ defined by $\varphi^{\mu}(\la)=\varphi\lr{a_{\la}a_{\mu}}$ is a finite harmonic function on $\Ga$, while $\varphi(\mu)<+\Inf$. Furthermore, the following inequality holds $\varphi\geq \const\cdot \varphi^{\mu}$. Next, observe that $\varphi^{\mu}$ is strictly positive, since $a_{\la}a_{\mu}\neq0$ and structure constants of $A$ are non-negative with respect to the basis $\{a_{\la}\}_{\la\in\Ga}$. Therefore, $\varphi$ and $\varphi^{\mu}$ are proportional. Thus, $\varphi$ is finite.
\end{proof}

\begin{corollary}\cite[p. 371, the paragraph just before Theorem 3.5]{boyer_symplectic}
If $\Ga$ admits a strictly positive indecomposable finite harmonic function, then it possesses no strictly positive semifinite indecomposable harmonic functions. 
\end{corollary}
\begin{proof}
Suppose that $\varphi$ is a strictly positive indecomposable finite harmonic function and $a_{\la}a_{\mu}=0$ for some $\la,\mu\in\Ga$. Then $\varphi(a_{\la}a_{\mu})=\varphi(0)=0$ and Theorem \ref{vershik kerov ring theorem} yields $\varphi(\la)\varphi(\mu)=0$, which contradicts the strict positivity of $\varphi$.
\end{proof}

\section{Boyer's Lemma}\label{5}
In this section we discuss a very useful claim related to arbitrary harmonic functions on a graded graph. It allows one to determine the finiteness ideal of an indecomposable semifinite harmonic function in several concrete situations. This principle, which was first observed by R. P. Boyer and published only in 1983, see \cite[Theorem 1.10, Example p.212]{boyer}, had been also stated by Wassermann \cite[Boyer's Lemma p.149]{wassermann1981} two years before the paper \cite{boyer}. We formulate and prove a slightly involved generalization of Wassermann's concise argument. It turns out to be a combinatorial analog of \cite[Theorem 1.10]{boyer}. After that we consider a couple of examples, which immediately follow from the general claim. Boyer's Lemma from \cite{wassermann1981} becomes a part of the first example, see Remark \ref{rem boyers classic}.

\subsection{General statement}Recall that the set of vertices lying on the $n$-th level of a graded graph $\Ga$ is denoted by $\Ga_n$. Below we work with arbitrary harmonic functions and do \textit{not} assume that they are finite or semifinite.

\begin{definition}\label{def semifin at a point}
A harmonic function $\varphi$ is called \textit{semifinite at a vertex} $\la$, if $\varphi(\la)=+\Inf$ and there exists a sequence $\{a_n\}_{n\geq 1}\subset \KK(\Ga)$ such that
\begin{itemize}
     \item $a_n\leq_K \la$,
     
     \item $\varphi(a_n)<+\Inf$,
     
     \item $\lim\limits_{n\to+\Inf}\varphi(a_n)=+\Inf$.
 \end{itemize} 
The sequence $\{a_n\}_{n\geq 1}$ will be called \textit{an approximating sequence for the vertex} $\la$.
\end{definition}

\begin{observation}\label{observation locally semifin}
If $\varphi$ is semifinite at a vertex $\la$, then for any vertex $\mu\leq\la$ the function $\varphi$ is semifinite at the vertex $\mu$ too.
\end{observation}

\begin{proposition}[Generalized Boyer's lemma]\label{gen boyer's lemma}
Let $(\Ga,\ka)$ be a graded graph and $\varphi$ be a harmonic function on it. Assume that $I\subset \Ga$ is an ideal, $J=\Ga\backslash I$ is the corresponding coideal and we are given a fixed vertex $\la\in J_n$. Suppose that there exists a natural number $m=m(\la)$ and a tuple of non-negative real numbers $\{\beta_{\nu}\}_{\nu\in I_m}$, which may depend on $\la$, such that the following conditions are satisfied
\begin{itemize}
    \item there exists a vertex $\nu\in I_m$ with $\beta_{\nu}\neq 0$ and $\varphi(\nu)>0$,
    
    \item for any large enough $l$ and any vertex $\eta\in I_{n+l+1}$ the following inequality holds
    \begin{speqn}\label{ineq gen boyers}
    \summ_{\mu\in J_{n+l}}\dim(\la,\mu)\ka(\mu,\eta)\geq{\summ_{\nu\in I_m}\beta_{\nu}\dim(\nu,\eta)}.
    \end{speqn}
\end{itemize}

Then $\varphi(\la)=+\Inf$. If in addition $\varphi(\nu)<+\infty$ for any $\nu\in I_m$ such that $\beta_{\nu}\neq 0$, then $\varphi$ is semifinite at the vertex $\la$.
\end{proposition}

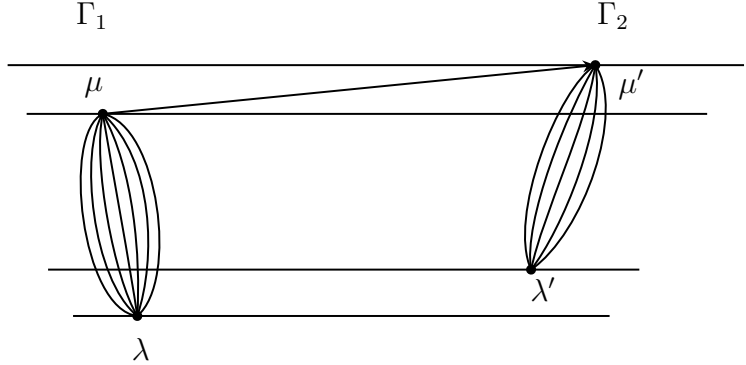
\begin{figure}[h]
    \centering
    \tikzset{every picture/.style={line width=0.75pt}} 

\begin{tikzpicture}[x=0.75pt,y=0.75pt,yscale=-0.6,xscale=0.6]

\draw    (129,130) .. controls (109,160.5) and (120,269.5) .. (158,300) ;
\draw [shift={(158,300)}, rotate = 38.75] [color={rgb, 255:red, 0; green, 0; blue, 0 }  ][fill={rgb, 255:red, 0; green, 0; blue, 0 }  ][line width=0.75]      (0, 0) circle [x radius= 3.35, y radius= 3.35]   ;
\draw [shift={(129,130)}, rotate = 123.25] [color={rgb, 255:red, 0; green, 0; blue, 0 }  ][fill={rgb, 255:red, 0; green, 0; blue, 0 }  ][line width=0.75]      (0, 0) circle [x radius= 3.35, y radius= 3.35]   ;
\draw    (65,130) -- (637,130) ;
\draw    (49,89) -- (668,89) ;
\draw    (104,300) -- (555,300) ;
\draw    (83,261) -- (580,261) ;
\draw    (129,130) .. controls (120,170.5) and (139,269.5) .. (158,300) ;
\draw    (129,130) .. controls (137,177.5) and (152,258.5) .. (158,300) ;
\draw    (129,130) .. controls (155,176.5) and (161,250.5) .. (158,300) ;
\draw    (129,130) .. controls (172,168.5) and (174,252.5) .. (158,300) ;
\draw    (543,89) .. controls (512,108.5) and (469,216.5) .. (489,261) ;
\draw [shift={(489,261)}, rotate = 65.8] [color={rgb, 255:red, 0; green, 0; blue, 0 }  ][fill={rgb, 255:red, 0; green, 0; blue, 0 }  ][line width=0.75]      (0, 0) circle [x radius= 3.35, y radius= 3.35]   ;
\draw [shift={(543,89)}, rotate = 147.83] [color={rgb, 255:red, 0; green, 0; blue, 0 }  ][fill={rgb, 255:red, 0; green, 0; blue, 0 }  ][line width=0.75]      (0, 0) circle [x radius= 3.35, y radius= 3.35]   ;
\draw    (543,89) .. controls (520,126.5) and (482,220.5) .. (489,261) ;
\draw    (543,89) .. controls (536,139.5) and (499,213.5) .. (489,261) ;
\draw    (543,89) .. controls (552,138.5) and (515,227.5) .. (489,261) ;
\draw    (543,89) .. controls (569,128.5) and (533,230.5) .. (489,261) ;
\draw    (129,130) .. controls (94,146.5) and (111,294.5) .. (158,300) ;
\draw    (129,130) .. controls (176,139.5) and (193,270.5) .. (158,300) ;
\draw    (540.01,89.3) -- (129,130) ;
\draw [shift={(543,89)}, rotate = 174.34] [fill={rgb, 255:red, 0; green, 0; blue, 0 }  ][line width=0.08]  [draw opacity=0] (10.72,-5.15) -- (0,0) -- (10.72,5.15) -- (7.12,0) -- cycle    ;

\draw (104,34) node [anchor=north west][inner sep=0.75pt]    {$\Gamma_{1}$};
\draw (542,34) node [anchor=north west][inner sep=0.75pt]    {$\Gamma_{2}$};
\draw (151,316) node [anchor=north west][inner sep=0.75pt]    {$\lambda$};
\draw (486,265) node [anchor=north west][inner sep=0.75pt]    {$\lambda'$};
\draw (111,97) node [anchor=north west][inner sep=0.75pt]    {$\mu$};
\draw (560,90) node [anchor=north west][inner sep=0.75pt]    {$\mu'$};

\end{tikzpicture}
    \caption{Condition \eqref{boyer's condition2} means that the "number" of paths from $\la$ to $\mu'$, that go through $\mu$, is not smaller then the "number" of arbitrary paths from $\la'$ to $\mu'$.}
    \label{picture example1}
    \end{figure}

\begin{remark}
Condition \eqref{ineq gen boyers} is a refinement of some condition on the "number" of paths in the graph $\Ga$, which admits a graphical interpretation, see condition \eqref{boyer's condition2} from Corollary \ref{cor boyers lemma} and Figure \ref{picture example1}.
\end{remark}
\begin{proof}[Proof of Proposition \ref{gen boyer's lemma}]
Let us multiply \eqref{ineq gen boyers} by $\eta\in K_0(\Ga)$ and sum over all $\eta\in I_{n+l+1}$. Then we get
\begin{speqn}\label{ineq gen boyers proof}
\summ_{\substack{\eta\in I_{n+l+1}\\ \mu\in J_{n+l}}}\dim(\la,\mu)\ka(\mu,\eta)\cdot\eta\geq_{K}\summ_{\substack{\eta\in I_{n+l+1} \\ \nu\in I_m}}\beta_{\nu}\dim(\nu,\eta)\cdot\eta,
\end{speqn}
where the both sides of the inequality are considered as elements of $\K(\Ga)$ and the partial order on $\K(\Ga)$ defined by the cone $\KK(\Ga)$ is denoted by $\geq_K$. Furthermore, the right-hand side of the inequality \eqref{ineq gen boyers proof} equals $\summ_{\nu\in I_m}\beta_{\nu}\nu$. Let us denote it by $b_{\la}$. Then $\varphi(b_{\la})>0$ and
$$\summ_{\substack{\eta\in I_{n+l+1}\\ \mu\in J_{n+l}}}\dim(\la,\mu)\ka(\mu,\eta)\eta\geq_{K}b_{\la}.$$ 

The only thing we are left to do is to reproduce the original argument of A. Wassermann \cite[p.149, the proof of Boyer's Lemma]{wassermann1981} in our context:

\begin{speqn}\label{ineq1}
    \la=\summ_{\ol{\eta}\in \Ga_{n+N+1}}\dim(\la,\ol{\eta})\ol{\eta}\geq_K \summ_{\ol{\eta}\in I_{n+N+1}}\dim(\la,\ol{\eta})\ol{\eta}.
\end{speqn}

Note that, if $\la\in J_n$ and $\ol{\eta}\in I_{n+N+1}$, then
\begin{speqn}\label{ineq2}
    \dim(\la,\ol{\eta})=\summ_{l=0}^{N}\summ_{\substack{\eta\in I_{n+l+1}\\ \mu\in J_{n+l}}}\dim(\la,\mu)\ka(\mu,\eta)\dim(\eta,\ol{\eta}).
\end{speqn}

Substitute \eqref{ineq2} into \eqref{ineq1}:     $\la\geq_K \summ_{l=0}^{N}\summ_{\substack{\eta\in I_{n+l+1}\\ \mu\in J_{n+l}}}\summ_{\ol{\eta}\in I_{n+N+1}}\dim(\la,\mu)\ka(\mu,\eta)\dim(\eta,\ol{\eta})\ol{\eta}$.

Now sum over $\ol{\eta}$: \begin{speqn}\label{askudfgyil}
\la\geq_K \summ_{l=0}^{N}\summ_{\substack{\eta\in I_{n+l+1}\\ \mu\in J_{n+l}}}\dim(\la,\mu)\ka(\mu,\eta)\eta\geq_K b_{\la}\cdot N
\end{speqn}

Compare \eqref{askudfgyil} with $(1.10.1)$ and $(1.10.2)$ from \cite[Theorem 1.10]{boyer}.

Thus, \eqref{askudfgyil} yields $\varphi(\la)\geq \varphi(b_{\la})\cdot N$ for any $N$ hence $\varphi(\la)=+\infty$. Moreover, the sequence $a_N=b_{\la}\cdot N$ is an approximating sequence for the vertex $\la$ if $\varphi(b_{\la})<+\Inf$.
\end{proof}

\subsection{Example 1}
Consider graded graphs $(\Ga_1,\ka_1)$ and $(\Ga_2,\ka_2)$ and suppose that we are given a graded map $\Ga_1\rightarrow \Ga_2$, $\la\mapsto \la'$. Let $\lr{\Ga,\ka}$ be one more graded graph that satisfies the following requirements:
\mathtoolsset{showonlyrefs=false}
\begin{equation}\label{c}
    \lr{\Ga}_n=\lr{\Ga_1}_n\sqcup \lr{\Ga_2}_{n-1}\ \text{for}\ n\geq 1,\ \lr{\Ga}_0=\lr{\Ga_1}_0.
\end{equation}

\begin{speqn}\label{2c}
    \ka\lr{\la,\mu}=\ka_1\lr{\la,\mu},\ \text{if}\ \la,\mu\in\Ga_1,\\
    \ka\lr{\la,\mu}=\ka_2\lr{\la,\mu},\ \text{if}\ \la,\mu\in\Ga_2.
\end{speqn}
\begin{speqn}\label{3c}
    \ka\lr{\la,\mu}=0,\ \text{if}\ \la\in\Ga_2,\mu\in\Ga_1.\\
\end{speqn}
\mathtoolsset{showonlyrefs=true}

\begin{figure}[h]
    \centering
    \tikzset{every picture/.style={line width=0.75pt}} 

\begin{tikzpicture}[x=0.75pt,y=0.75pt,yscale=-0.7,xscale=0.7]

\draw    (157,207.5) -- (407,207.5) ;
\draw    (103,118) -- (550,118) ;
\draw    (145.82,120.38) -- (212.5,207.5) ;
\draw [shift={(212.5,207.5)}, rotate = 52.57] [color={rgb, 255:red, 0; green, 0; blue, 0 }  ][fill={rgb, 255:red, 0; green, 0; blue, 0 }  ][line width=0.75]      (0, 0) circle [x radius= 3.35, y radius= 3.35]   ;
\draw [shift={(144,118)}, rotate = 52.57] [fill={rgb, 255:red, 0; green, 0; blue, 0 }  ][line width=0.08]  [draw opacity=0] (10.72,-5.15) -- (0,0) -- (10.72,5.15) -- (7.12,0) -- cycle    ;
\draw    (208.15,121) -- (212.5,207.5) ;
\draw [shift={(208,118)}, rotate = 87.12] [fill={rgb, 255:red, 0; green, 0; blue, 0 }  ][line width=0.08]  [draw opacity=0] (10.72,-5.15) -- (0,0) -- (10.72,5.15) -- (7.12,0) -- cycle    ;
\draw    (269.36,120.51) -- (212.5,207.5) ;
\draw [shift={(271,118)}, rotate = 123.17] [fill={rgb, 255:red, 0; green, 0; blue, 0 }  ][line width=0.08]  [draw opacity=0] (10.72,-5.15) -- (0,0) -- (10.72,5.15) -- (7.12,0) -- cycle    ;
\draw    (500.13,118.88) -- (212.5,207.5) ;
\draw [shift={(503,118)}, rotate = 162.88] [fill={rgb, 255:red, 0; green, 0; blue, 0 }  ][line width=0.08]  [draw opacity=0] (10.72,-5.15) -- (0,0) -- (10.72,5.15) -- (7.12,0) -- cycle    ;
\draw (205,221) node [anchor=north west][inner sep=0.75pt]    {$\lambda $};
\draw (490,135) node [anchor=north west][inner sep=0.75pt]    {$\lambda '$};
\draw (190,80) node [anchor=north west][inner sep=0.75pt]    {$\Gamma _{1}$};
\draw (490,80) node [anchor=north west][inner sep=0.75pt]    {$\Gamma _{2}$};

\end{tikzpicture}
    \caption{Example of the branching rule for $\Ga$.}
    \label{branching rule for Boyers classic}
\end{figure}
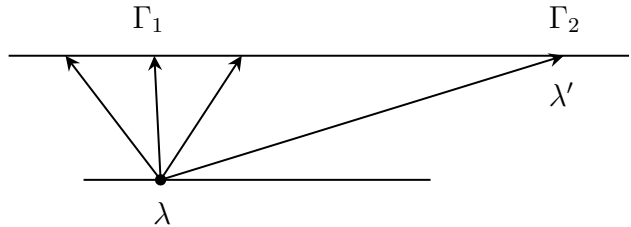

Condition \eqref{3c} means that $\Ga_2$ is an ideal of $\Ga$. For simplicity one can assume that edges from $\Ga_1$ to $\Ga_2$ can go from $\la$ to $\la'$ only, see Figure \ref{branching rule for Boyers classic}. But we will not use this later on.

    \begin{corollary}\label{cor boyers lemma}
    Assume that the map $\nu\mapsto\nu'$ is surjective and let $\la\in \lr{\Ga_1}_n$ be a fixed vertex. Suppose that for any large enough $l$ and any vertex $\mu\in\lr{\Ga_1}_{n+l}$ the following inequality holds
    \begin{speqn}\label{boyer's condition2}
    \dim_1(\la,\mu)\ka(\mu,\mu')\geq \dim_2(\la',\mu'),
    \end{speqn}
    where $\dim_1(\cdot,\cdot)$ and $\dim_2(\cdot,\cdot)$ are shifted dimensions for $(\Ga_1,\ka_1)$ and $(\Ga_2,\ka_2)$. Now let $\varphi$ be a harmonic function on $\Ga$ with $\varphi(\la')>0$. Then $\varphi(\la)=+\infty$ and $\varphi$ is semifinite at the vertex $\la$ if $\varphi(\la')<+\Inf$.
    \end{corollary}
    \begin{proof}
    Recall that $\Ga_2$ is an ideal of $\Ga$. Therefore, we may apply Proposition \ref{gen boyer's lemma} for $I=\Ga_2$, $J=\Ga_1$, $m=|\la|+1$ and $\beta_{\nu}=\delta_{\nu,\la'}$. Then we bound from below the sum in the left hand side of \eqref{ineq gen boyers} in terms of one of its summands and use \eqref{boyer's condition2}.
    \end{proof}
    
    \begin{remark}\label{rem boyers classic}
    If the map $\la\mapsto \la'$ is a branching graph morphism, that is $\ka(\la,\mu)=\ka(\la',\mu')$, then condition \eqref{boyer's condition2} means that $\ka(\mu,\mu')\geq 1$. If the equality holds identically, then we obtain the original formulation of Boyer's Lemma \cites[p.149, Boyer's lemma]{wassermann1981}.
    \end{remark}
    
\subsection{Example 2}  
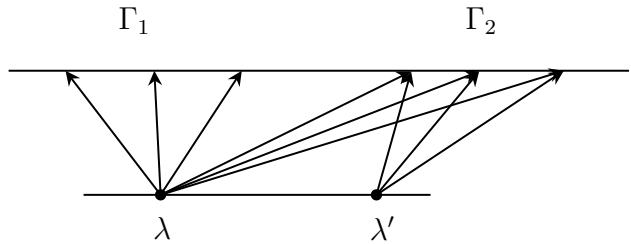
\begin{figure}[h]
    \centering
    \tikzset{every picture/.style={line width=0.75pt}} 

\begin{tikzpicture}[x=0.75pt,y=0.75pt,yscale=-0.7,xscale=0.7]

\draw    (177,227.5) -- (427,227.5) ;
\draw    (123,138) -- (570,138) ;
\draw    (165.82,140.38) -- (232.5,227.5) ;
\draw [shift={(232.5,227.5)}, rotate = 52.57] [color={rgb, 255:red, 0; green, 0; blue, 0 }  ][fill={rgb, 255:red, 0; green, 0; blue, 0 }  ][line width=0.75]      (0, 0) circle [x radius= 3.35, y radius= 3.35]   ;
\draw [shift={(164,138)}, rotate = 52.57] [fill={rgb, 255:red, 0; green, 0; blue, 0 }  ][line width=0.08]  [draw opacity=0] (10.72,-5.15) -- (0,0) -- (10.72,5.15) -- (7.12,0) -- cycle    ;
\draw    (228.15,141) -- (232.5,227.5) ;
\draw [shift={(228,138)}, rotate = 87.12] [fill={rgb, 255:red, 0; green, 0; blue, 0 }  ][line width=0.08]  [draw opacity=0] (10.72,-5.15) -- (0,0) -- (10.72,5.15) -- (7.12,0) -- cycle    ;
\draw    (289.36,140.51) -- (232.5,227.5) ;
\draw [shift={(291,138)}, rotate = 123.17] [fill={rgb, 255:red, 0; green, 0; blue, 0 }  ][line width=0.08]  [draw opacity=0] (10.72,-5.15) -- (0,0) -- (10.72,5.15) -- (7.12,0) -- cycle    ;
\draw    (520.13,138.88) -- (232.5,227.5) ;
\draw [shift={(523,138)}, rotate = 162.88] [fill={rgb, 255:red, 0; green, 0; blue, 0 }  ][line width=0.08]  [draw opacity=0] (10.72,-5.15) -- (0,0) -- (10.72,5.15) -- (7.12,0) -- cycle    ;
\draw    (520.5,139.66) -- (388,227.5) ;
\draw [shift={(388,227.5)}, rotate = 146.46] [color={rgb, 255:red, 0; green, 0; blue, 0 }  ][fill={rgb, 255:red, 0; green, 0; blue, 0 }  ][line width=0.75]      (0, 0) circle [x radius= 3.35, y radius= 3.35]   ;
\draw [shift={(523,138)}, rotate = 146.46] [fill={rgb, 255:red, 0; green, 0; blue, 0 }  ][line width=0.08]  [draw opacity=0] (10.72,-5.15) -- (0,0) -- (10.72,5.15) -- (7.12,0) -- cycle    ;
\draw    (412.19,141.39) -- (388,227.5) ;
\draw [shift={(413,138.5)}, rotate = 105.69] [fill={rgb, 255:red, 0; green, 0; blue, 0 }  ][line width=0.08]  [draw opacity=0] (10.72,-5.15) -- (0,0) -- (10.72,5.15) -- (7.12,0) -- cycle    ;
\draw    (410.31,139.83) -- (232.5,227.5) ;
\draw [shift={(413,138.5)}, rotate = 153.75] [fill={rgb, 255:red, 0; green, 0; blue, 0 }  ][line width=0.08]  [draw opacity=0] (10.72,-5.15) -- (0,0) -- (10.72,5.15) -- (7.12,0) -- cycle    ;
\draw    (460.08,140.81) -- (388,227.5) ;
\draw [shift={(462,138.5)}, rotate = 129.74] [fill={rgb, 255:red, 0; green, 0; blue, 0 }  ][line width=0.08]  [draw opacity=0] (10.72,-5.15) -- (0,0) -- (10.72,5.15) -- (7.12,0) -- cycle    ;
\draw    (459.2,139.58) -- (232.5,227.5) ;
\draw [shift={(462,138.5)}, rotate = 158.8] [fill={rgb, 255:red, 0; green, 0; blue, 0 }  ][line width=0.08]  [draw opacity=0] (10.72,-5.15) -- (0,0) -- (10.72,5.15) -- (7.12,0) -- cycle    ;

\draw (225,241) node [anchor=north west][inner sep=0.75pt]    {$\lambda $};
\draw (381,241) node [anchor=north west][inner sep=0.75pt]    {$\lambda '$};
\draw (200,90) node [anchor=north west][inner sep=0.75pt]    {$\Gamma _{1}$};
\draw (450,90) node [anchor=north west][inner sep=0.75pt]    {$\Gamma _{2}$};

\end{tikzpicture}
    \caption{Example of the branching rule for $\Ga$.}
    \label{branching rule for twisted Boyers classic}
\end{figure}

Let us consider graded graphs $(\Ga_1,\ka_1)$ and $(\Ga_2,\ka_2)$ and suppose that we are given a graded map $\Ga_1\rightarrow \Ga_2$, $\la\mapsto \la'$. Let $\lr{\Ga,\ka}$ be another graded graph, that satisfies $\lr{\Ga}_n=\lr{\Ga_1}_n\sqcup \lr{\Ga_2}_{n}\ \text{for}\ n\geq 0$, and conditions \eqref{2c}, \eqref{3c}. Recall that the last condition means that $\Ga_2$ is an ideal of $\Ga$. For simplicity one can assume that vertices $\la\in\Ga_1$ and $\mu\in\Ga_2$ are joined by an edge if and only if  $\la'\nearrow \mu$, as it is shown on Figure \ref{branching rule for twisted Boyers classic}. 

\begin{corollary}
Suppose that the map $\la\mapsto\la'$ is surjective. Let $\la\in\Ga_1$ be a fixed vertex and assume that the following inequalities hold for any $\mu\in\Ga_1$
\begin{speqn}
    \ka(\la,\mu)\geq\ka(\la',\mu'),\\
    \ka(\la,\mu')\geq \ka(\la',\mu').
\end{speqn}
Then $\varphi(\la)=+\Inf$ for any harmonic function $\varphi$ on $\Ga$ such that $\varphi(\la')>0$. Moreover, $\varphi$ is semifinite at the vertex $\la$, if $0<\varphi(\la')<+\Inf$.
\end{corollary}
\begin{proof}
Let us take $I=\Ga_2$, $J=\Ga_1$, $m=|\la|$ and $\beta_{\nu}=\delta_{\nu,\la'}$ in Proposition \ref{gen boyer's lemma} and prove that the following inequality holds $\summ_{\mu\in\Ga_1}\dim(\la,\mu)\ka(\mu,\eta)\geq\dim(\la',\eta)$ for any $\eta\in\Ga_2$. In order to do so, we check that $\dim(\la,\mu)\geq \dim(\la',\mu')$ and write
\begin{speqn}
\cfrac{\summ_{\mu\in\Ga_1}\dim(\la,\mu)\ka(\mu,\eta)}{\dim(\la',\eta)}\geq \cfrac{\summ_{\mu\in\Ga_1}\dim(\la',\mu')\ka(\mu',\eta)}{\dim(\la',\eta)}\geq \cfrac{\summ_{\ol{\mu}\in\Ga_2}\dim(\la',\ol{\mu})\ka(\ol{\mu},\eta)}{\dim(\la',\eta)}=1.
\end{speqn}
For each of these inequalities we have used that $\la\mapsto \la'$ is surjective.
\end{proof}

\begin{remark}\label{rem mac}
As it was pointed out in the introduction, one can obtain an exhaustive list of indecomposable semifinite harmonic functions on the Macdonald graph, which corresponds to the simplest Pieri rule for the Macdonald symmetric functions, by applying Wassermann's method. This list turns out to be very similar to that for the Young graph, see \cite[Theorem 9 on page 150]{wassermann1981}. For instance, the space of classification parameters is an obvious $(q,t)$-deformation of the parameter space for the Young graph. Namely, we should deform only the continuous part of the data in the same way as it is deformed in the case of finite harmonic functions, replacing the ordinary Thoma simplex with the $(q,t)$-deformed Thoma simplex, see Theorem 1.4 and Proposition 1.6 from \cite{matveev2018macdonaldpositive}, while the discrete part remains the same. This result easily follows from the original argument of A. Wassermann, Theorem 1.4 and Proposition 1.6 from \cite{matveev2018macdonaldpositive}, Proposition \ref{co-ideal}, Theorem \ref{theorem ext}, Proposition \ref{prop prim coideal}, and Corollary \ref{cor boyers lemma}. Instead of using Theorem \ref{theorem mult graph} we must apply a similar argument obtained with the help of a trick due to K. Matveev \cite[§6,  Proof of Proposition 1.6]{matveev2018macdonaldpositive}.
\end{remark}

\appendix
\numberwithin{equation}{section}

\section{Direct product of branching graphs}\label{addendum}

In this appendix we describe indecomposable finite harmonic functions on the product of branching graphs in terms of harmonic functions on the multipliers. This result is not related to semifinite harmonic functions in a straightforward way, but it turns out to be very useful for the description of semifinite harmonic functions on some branching graphs such as the Gnedin-Kingman graph \cite{safonkin20} and the zigzag graph. The latter was studied in the paper \cite{gnedin_olsh2006}. One can treat the main result of this appendix, Proposition \ref{product graph}, as a generalization of the well known de Finetti theorem \cite[Theorem 5.1, Theorem 5.2]{bor_olsh2017}. The difference between Proposition \ref{product graph} ($n=2$ case) and the de Finetti theorem is that we replace two sides of the Pascal triangle, which correspond to two embeddings $\Z\hookrightarrow \Z\oplus \Z$ along the first and the second components, with arbitrary branching graphs. Remark that the case when one of these graphs is a line consisting of one vertex at each level has been already known, see \cite[Theorem 2.8]{versh_nikitin2010}. Note that in this theorem one should consider only strictly positive harmonic functions (or, equivalently, central measures) instead of arbitrary ones.

Let us provide some motivation for the main definition of the present section. If $A$ and $B$ are unital $\ZZ$-graded $\R$-algebras, then their tensor product (over $\R$) is a unital graded algebra too. Namely, if $A=\bigoplus\limits_{n\geq 0}A_n$, $A_0=\R$ and $B=\bigoplus\limits_{n\geq 0}B_n$, $B_0=\R$, then $A\otimes_{\R} B=\bigoplus\limits_{k\geq 0}\lr{A\otimes_{\R}B}_k$, where
$$\lr{A\otimes_{\R}B}_k=\bigoplus\limits_{\substack{n,m\geq 0:\\ n+m=k}}A_n\otimes_{\R} B_m.$$ 
Furthermore, $\mathds{1}_{A\otimes B}=\mathds{1}_A\otimes \mathds{1}_B$ and $\lr{A\otimes_{\R}B}_0=\R\cdot \mathds{1}_{A\otimes B}$. This simple fact, together with Definition \ref{multiplicative def}, motivates us to consider the \textit{direct product} of two graded graphs.
 	
\begin{definition}\label{direct product graph def}
By the \textit{direct product} of graded graphs $\lr{\Ga_1,\ka_1}$ and $\lr{\Ga_2,\ka_2}$ we mean the graded graph $\lr{\Ga_1\times \Ga_2,\ka_{1}\times \ka_{2}}$, where
$$\lr{\Ga_1\times\Ga_1}_k=\bigsqcup\limits_{\substack{n,m\geq0: \\ n+m=k}}\lr{\Ga_1}_n\times\lr{\Ga_2}_m$$
and 
$$\lr{\ka_1\times\ka_2}\biglr{(\la_1,\mu_1);(\la_2,\mu_2)}=
\begin{cases}
    \ka_1(\la_1,\la_2),\ &\text{if}\ \mu_1=\mu_2,\\
    \ka_2(\mu_1,\mu_2),\ &\text{if}\ \la_1=\la_2,\\
    0\ &\text{otherwise.}
\end{cases}$$
\end{definition} 
	
The next lemma ties together some properties of the direct product of graded graphs. 
    
The subset $\Ga_{\la}=\{\mu\in\Ga\mid \mu\leq\la\}$ of a graded graph $\Ga$ is called \textit{the principle coideal} associated to $\la\in\Ga$.
	
\begin{lemma}\label{prop prod}
Let $\Ga_1$ and $\Ga_2$ be graded graphs. 
\begin{enumerate}[leftmargin=5ex,label=\theenumi)]
	\item A graph $\Ga_1\times\Ga_2$ is primitive if and only if  $\Ga_1$ and $\Ga_2$ are primitive.
	\vspace{0.5\baselineskip}
	\item If $\Ga_1$ and $\Ga_2$ are branching graphs and $J\subset\Ga_1\times\Ga_2$ is a saturated primitive coideal, then there exist coideals $J_1\subset\Ga_1$ and $J_2\subset \Ga_2$ such that $J=J_1\times J_2$ and 
	\begin{itemize}
	\item $J_1,J_2$ are saturated and primitive or
	        
	\item $J_1$ is principle and $J_2$ is saturated and primitive or
	        
	\item $J_1$ is saturated and primitive and $J_2$ is principle. 
	\end{itemize}
	Moreover, coideals $J_1$ and $J_2$ are uniquely defined. 
	\vspace{0.5\baselineskip}
	\item Let $\la,\la'\in\Ga_1$ and $\mu,\mu'\in\Ga_2$. Then
	\begin{speqn}
	\dim\biglr{(\la,\mu),(\la',\mu')}=\binom{|\la'|-|\la|+|\mu'|-|\mu|}{|\la'|-|\la|} \dim_1\lr{\la,\la'}\dim_2\lr{\mu,\mu'},
	\end{speqn}
	where $\binom{n}{k}$ denotes the binomial coefficient and $\dim_1(\cdot,\cdot)$, $\dim_2(\cdot,\cdot)$ are shifted dimensions for $\Ga_1$ and $\Ga_2$, see \eqref{shiftdim} on page \pageref{shiftdim}.
	\end{enumerate}
	\end{lemma}
	
\begin{proof} The first and the second assertions follow from Proposition \ref{co-ideal} immediately and the third one is obvious.
\end{proof}

Note that we can easily generalize the statement of Lemma \ref{prop prod} to the case of $n>2$ graded graphs. Furthermore, the direct product of multiplicative graphs is multiplicative too. For the direct product of two multiplicative graphs the corresponding algebra is the tensor product of the initial algebras, the distinguished basis is the tensor product of the bases and the element that was denoted by $\wh{a}$ in Definition \ref{multiplicative def} is $\wh{a}\otimes_{\R}\mathds{1}_B+\mathds{1}_A\otimes_{\R} \wh{b}$, where $\wh{a}$ and $\wh{b}$ are the same elements for the initial algebras. Thus, we can define the direct product of finitely many graded graphs and the product of multiplicative graphs is multiplicative as well.

Recall that a harmonic function $\varphi$ on a branching graph $\Ga$ is called normalized if $\varphi(\diameter)=1$. 

\begin{remark}
Let $\Ga_1,\ldots,\Ga_n$ be branching graphs and let $\varphi_1,\ldots,\varphi_n$ be finite normalized harmonic functions on them. Then the function $\varphi\colon\Ga_1\times\ldots\times\Ga_n\rightarrow \R_{\geq0}$ defined by
\begin{speqn}\label{injection of harm funct}
\varphi(\la_1,\ldots,\la_n)=w_1^{|\la_1|}\ldots w_n^{|\la_n|}\varphi_1(\la_1)\ldots\varphi_n(\la_n)
\end{speqn}
is harmonic and normalized whenever $w_1,\ldots,w_n\in\RR$ and $w_1+\ldots+w_n=1$. 

Remark that we can recover these $\varphi_1,\ldots, \varphi_n$ and $w_1,\ldots, w_n$ from $\varphi$ as follows. Let us set $$\MC(a_1,\ldots,a_n)=\binom{a_1+\ldots+a_n}{a_1,\ldots,a_n}=\cfrac{(a_1+\ldots+a_n)!}{a_1!\ldots a_n!}.$$

Then
\begin{equation}\label{as8}
    \begin{multlined}
    \varphi_i(\mu)=\summ_{\substack{\la_j\in\Ga_j, j\neq i\\ j=1,\ldots,n}}\MC\biglr{|\la_1|,\ldots,|\la_{i-1}|,|\mu|-1,|\la_{i+1}|,\ldots,|\la_n|}\ \ \ \ \ \ \ \ \ \ \ \ \ \ \ \ \ \ \ \ \ \ \\
     \cdot\prodd_{\substack{j=1\\ j\neq i}}^{n}\dim(\la_j)\cdot  \varphi(\la_1,\ldots,\la_{i-1},\mu,\la_{i+1},\ldots,\la_n)
    \end{multlined}
\end{equation}
for $|\mu|\geq 1$ and 
\begin{equation}\label{as9}
    w_1^{k_1}\ldots w_n^{k_n}=\summ_{\substack{\la_i\in\Ga_i, |\la_i|=k_i\\ i=1,\ldots,n}}\dim(\la_1)\ldots\dim(\la_n)\cdot\varphi(\la_1,\ldots,\la_n)
\end{equation}
for any positive integers $k_1,\ldots,k_n$. 

Compare \eqref{as8} and \eqref{as9} with the first two formulas from the proof of Theorem 2.8 in \cite{versh_nikitin2010}.

\end{remark}

\begin{notation}
Let $\lr{\Ga,\ka}$ be a branching graph. We denote by $\exFH(\Ga)$ the set of all finite normalized harmonic functions on $\Ga$ and by $\exFHo(\Ga)$  the subset of all strictly positive functions.
\end{notation}

\begin{proposition}\label{product graph}
Let $\Ga_1,\ldots, \Ga_n$ be branching graphs and $\Delta_n^{0}$ be the interior of the $n-1$-dimensional simplex, i.e. $\Delta_n^{0}=\{\lr{w_1,\ldots,w_n}\mid w_1+\ldots+w_n=1,\ w_i>0\}.$
\begin{enumerate}[leftmargin=3ex,label=\theenumi)]
    \item There is a bijection between $\exFHo(\Ga_1\times\ldots\times\Ga_n)$ and $\exFHo(\Ga_1)\times\ldots\times\exFHo(\Ga_n)\times \Delta_n^{0}$ defined by \eqref{injection of harm funct}.
    \vspace{1\baselineskip}
    \item There is a bijection between $\exFH(\Ga_1\times\ldots\times\Ga_n)$ and $\bigsqcup\limits_{\substack{I\colon I\subset \{1,2,\ldots, n\}\\ I\neq \diameter}}\Delta_{|I|}^{0}\times \bigtimes\limits_{i\in I}\exFH(\Ga_i)$.
    
    More precisely, for any harmonic function $\varphi\in\exFH(\Ga_1\times\ldots\times\Ga_n)$ there exist a non-empty set $I\subset \{1,2,\ldots,n\}$, harmonic functions $\varphi_i\in\exFH(\Ga_i)$, which are indexed by $i\in I$, and $w\in\Delta_{|I|}^{0}$ such that for any $n$-tuple of vertices $\la_1\in\Ga_1, \ldots, \la_n\in\Ga_n$ the following identity holds
    \begin{speqn}
    \varphi\lr{\la_1,\ldots,\la_n}=\begin{cases}
    \prodd_{i\in I}w_i^{|\la_i|}\varphi_i(\la_i),\ \text{if}\ \la_j=\diameter,\ \forall j\in\{1,2,\ldots,n\}\backslash I,\\
    0\ \text{otherwise}.
    \end{cases}
    \end{speqn}
    
    Moreover, these $I$, $\varphi_i$ and $w$ are uniquely defined.
    \end{enumerate}
    \end{proposition}
    \begin{remark}\label{remmn}
For multiplicative graphs Proposition \ref{product graph} is a straightforward consequence of the Vershik-Kerov ring theorem (Theorem \ref{vershik kerov ring theorem}). Namely, we should apply this theorem to the following simple fact
$$\Hom\biglr{A_1\otimes_{\R}\ldots\otimes_{\R} A_n,\R}\simeq \bigtimes\limits_{i=1}^{n}\Hom\lr{A_i,\R},$$ 
where $\Hom$ stands for the set of algebra homomorphisms. Indeed, to prove the first part of the proposition we note that there are two mutually inverse maps
\begin{speqn}
\Phi_{\raisebox{-0.5ex}{$\rightarrow$}}\colon\exFHo(\Ga_1)\times\ldots\times\exFHo(\Ga_n)\times \Delta_n^{0}\longrightarrow \exFHo(\Ga_1\times\ldots\times\Ga_n),\\
\lr{\varphi_1, \ldots,\varphi_n,w}\mapsto \lr{\varphi_1\circ r_{w_1}}\otimes \ldots \otimes \lr{\varphi_n\circ r_{w_n}};
\end{speqn}
and

\begin{speqn}
\Phi_{\raisebox{-0.5ex}{$\leftarrow$}}\colon\exFHo(\Ga_1\times\ldots\times\Ga_n)\longrightarrow\exFHo(\Ga_1)\times\ldots\times\exFHo(\Ga_n)\times \Delta_n^{0},\\
\varphi\mapsto\lr{\varphi\left.\right|_{A_1}\circ r^{-1}_{w_1},\ldots,\varphi\left.\right|_{A_n}\circ r^{-1}_{w_n},w}.
\end{speqn}
    
Here $r_u$ denotes the automorphism of a graded algebra defined on homogeneous elements as $a\mapsto u^{\deg{a}}a$, and $\varphi\left.\right|_{A_i}$ is the restriction of $\varphi\colon A_1\otimes\ldots\otimes A_n\rightarrow \R$ to the subalgebra $1^{\otimes i-1}\otimes A_i\otimes 1^{\otimes n-i}\simeq A_i$. Furthermore, the $n$-tuple $w=\lr{w_1,\ldots,w_n}$ that appears in the definition of the map $\Phi_{\raisebox{-0.5ex}{$\leftarrow$}}$ has the following form $w_i=\varphi\lr{1^{\otimes i-1}\otimes \wh{a}^{(i)}\otimes 1^{\otimes n-i}}$. Recall that the element $\wh{a}^{(i)}\in A_i$ defines the branching rule for $\Ga_i$, see Definition \ref{multiplicative def}.
\end{remark}

\begin{proof}[Proof of Proposition \ref{product graph}]
We prove the first part of the proposition for $n=2$ only. The case $n>2$ can be dealt with in the same manner. One can prove the second part of the proposition applying essentially the same argument, Proposition \ref{prop prim coideal} and the second part of Lemma \ref{prop prod}.
    
One can check that for any harmonic function $\varphi$ on $\Ga_1\times\Ga_2$ the right-hand side of \eqref{as8} defines a harmonic function on $\Ga_i$. Thus, $\varphi$ defined by \eqref{injection of harm funct} is indecomposable if $\varphi_1$ and $\varphi_2$ are. Then \eqref{injection of harm funct} defines an injective map $\exFHo(\Ga_1)\times\exFHo(\Ga_2)\times\Delta_2^{0}\longrightarrow \exFHo(\Ga_1\times\Ga_2)$. With the help of the Vershik-Kerov ergodic method, see \cites[p.20 Theorem 2]{versh_ker_81}[p.60, Theorem]{kerov_book}, we will show that this map is also surjective. Let $\varphi$ be a finite strictly positive normalized indecomposable harmonic function on $\Ga_1\times\Ga_2$. Then by \cite[p.60, Theorem]{kerov_book} there exists a path $\tau=\biglr{(\diameter,\diameter),(\la_1,\mu_1),\ldots}\in\T(\Ga_1\times\Ga_2)$ such that
$$\varphi(\la,\mu)=\lim\limits_{N\to+\Inf}\cfrac{\dim\biglr{(\la,\mu),(\la'_N,\mu'_N)}}{\dim\biglr{(\la'_N,\mu'_N)}}.$$ 
From the last part of Lemma \ref{prop prod} it follows that 
$$\cfrac{\dim\biglr{(\la,\mu),(\la'_N,\mu'_N)}}{\dim\biglr{(\la'_N,\mu'_N)}}=\cfrac{\biglr{|\la'_N|}^{\downarrow|\la|}\cdot\ \biglr{|\mu'_N|}^{\downarrow|\mu|}}{\biglr{|\la'_N|+|\mu'_N|}^{\downarrow (|\la|+|\mu|)}}\cdot\cfrac{\dim_1\lr{\la,\la'_N}}{\dim_1\lr{\la'_N}}\cdot\cfrac{\dim_2\lr{\mu,\mu'_N}}{\dim_2\lr{\mu'_N}},$$
where $x^{\downarrow k}=x(x-1)\ldots (x-k+1)$. Then strict positivity of $\varphi$ implies that $|\la'_N|\to+\Inf$ and $|\mu'_N|\to+\Inf$ as $N\to+\Inf$. Therefore, passing to appropriate subsequences we may assume that the following limits exist
$$\lim\limits_{N\to+\Inf}\cfrac{\dim_1\lr{\la,\la'_N}}{\dim_1\lr{\la'_N}},\ \lim\limits_{N\to+\Inf}\cfrac{\dim_2\lr{\mu,\mu'_N}}{\dim_2\lr{\mu'_N}},\ \lim\limits_{N\to+\Inf}\cfrac{|\la'_N|}{|\la'_N|+|\mu'_N|},\ \lim\limits_{N\to+\Inf}\cfrac{|\mu'_N|}{|\la'_N|+|\mu'_N|}.$$ 
Denoting them by $\varphi_1(\la)$, $\varphi_2(\mu)$, $w_1$ and $w_2$, we obtain the desired element of $\exFHo(\Ga_1)\times\exFHo(\Ga_2)\times\Delta_2^{0}$. Note that these $\varphi_1$ and $\varphi_2$ are indecomposable, since $\varphi$ is indecomposable.
\end{proof}
    
\begin{example}\label{example product graph}
Let us take $\Ga_1=\ldots=\Ga_n=\ZZ$ and assume that all edges are simple and go from $k$ to $k+1$ for $k\geq 0$. Then $\exFHo(\Ga_i)=\exFH(\Ga_i)$ is a singleton and $\Ga_1\times\ldots\times \Ga_n$ is the Pascal pyramid $\Pas_n$. Then from Proposition \ref{product graph} it follows that $\exFHo(\Pas_n)=\Delta_n^0$ and $\exFH(\Pas_n)= \bigsqcup\limits_{\substack{I\colon I\subset \{1,2,\ldots, n\}\\ I\neq \diameter}}\Delta_{|I|}^{0}=\Delta_n$, which is the $n-1$-dimensional simplex.
\end{example}
     
\begin{remark}
Proposition \ref{product graph} gives us the following view on Kerov's construction \cite[§4]{gnedin_olsh2006}. Comultiplication provides us a linear map $\K(\Ga)\rightarrow \K(\underbrace{\Ga\times\ldots\times\Ga}_n)$ and we take the composite of this map with an indecomposable harmonic function on $\underbrace{\Ga\times\ldots\times\Ga}_n$ to obtain an indecomposable harmonic function on $\Ga$.
\end{remark}

    \printbibliography

@book {bor_olsh2017,
    AUTHOR = {Borodin, Alexei and Olshanski, Grigori},
     TITLE = {Representations of the infinite symmetric group},
    SERIES = {Cambridge Studies in Advanced Mathematics},
    VOLUME = {160},
 PUBLISHER = {Cambridge University Press, Cambridge},
      YEAR = {2017},
     PAGES = {vii+160},
      %ISBN = {978-1-107-17555-6},
   MRCLASS = {20C32 (05E10 05E45 22D10 60J50 60K35)},
  MRNUMBER = {3618143},
MRREVIEWER = {Sevak Mkrtchyan},
}

@article {gnedin_olsh2006,
    AUTHOR = {Gnedin, Alexander and Olshanski, Grigori},
     TITLE = {Coherent permutations with descent statistic and the boundary
              problem for the graph of zigzag diagrams},
   JOURNAL = {Int. Math. Res. Not.},
  FJOURNAL = {International Mathematics Research Notices},
      YEAR = {2006},
     PAGES = {Art. ID 51968, 39},
      %ISSN = {1073-7928},
}

@article{kerov_vershik1990,
 title={The {G}rothendieck group of the infinite symmetric group and symmetric functions with the elements of the {$K_0$}-functor theory of {AF}-algebras},
  publisher = {Gordon and Breach},
  author={Kerov, S. and Vershik, A.},
  journal={Representation of Lie groups and related topics, Adv. Stud. Contemp. Math},
  volume={7},
  pages={36--114},
  year={1990}
}

@article{Bratteli1972,
 %ISSN = {00029947},
 abstract = {Inductive limits of ascending sequences of finite dimensional {$C^*$}-algebras are studied. The ideals of such algebras are classified, and a necessary and sufficient condition for isomorphism of two such algebras is obtained. The results of Powers concerning factor states and representations of UHF-algebras are generalized to this case. A study of the current algebra of the canonical anticommutation relations is then being made.},
 author = {Ola Bratteli},
 journal = {Transactions of the American Mathematical Society},
 pages = {195--234},
 publisher = {American Mathematical Society},
 title = {Inductive Limits of Finite Dimensional {$C^*$}-Algebras},
 volume = {171},
 year = {1972}
}

@phdthesis{wassermann1981,
author={Wassermann, Antony J.},
title={Automorphic actions of compact groups on operator algebras},
school={ University of Pennsylvania},
year={1981},
type = {{Ph.D.} thesis},
url={https://repository.upenn.edu/dissertations/AAI8127086/}
}

@book {strat_voic1975,
    AUTHOR = {Str\u{a}til\u{a}, \c{S}erban and Voiculescu, Dan},
     TITLE = {Representations of {AF}-algebras and of the group {$U(\infty
              )$}},
    SERIES = {Lecture Notes in Mathematics, Vol. 486},
 PUBLISHER = {Springer-Verlag, Berlin-New York},
      YEAR = {1975},
     PAGES = {viii+169},
}

@article {vershik_Kerov83,
    AUTHOR = {Kerov, S. V. and Vershik, A. M.},
     TITLE = {The {$K$}-functor ({G}rothendieck group) of the infinite
              symmetric group},
      %NOTE = {Differential geometry, Lie groups and mechanics, V},
   JOURNAL = {Zap. Nauchn. Sem. Leningrad. Otdel. Mat. Inst. Steklov.
              (LOMI)},
  FJOURNAL = {Zapiski Nauchnykh Seminarov Leningradskogo Otdeleniya
              Matematicheskogo Instituta imeni V. A. Steklova Akademii Nauk
              SSSR (LOMI)},
    VOLUME = {123},
      YEAR = {1983},
     PAGES = {126--151},
      %ISSN = {0373-2703},
   %MRCLASS = {20C32 (18F25 18F30 46L99 46M20)},
  %MRNUMBER = {697247},
%MRREVIEWER = {L\'{a}szl\'{o} Zsid\'{o}},
}

@book {kerov_book,
    AUTHOR = {Kerov, S. V.},
     TITLE = {Asymptotic representation theory of the symmetric group and
              its applications in analysis},
    SERIES = {Translations of Mathematical Monographs},
    VOLUME = {219},
    NOTE = {Translated from the Russian manuscript by N. V. Tsilevich,
              With a foreword by A. Vershik and comments by G. Olshanski},
 PUBLISHER = {American Mathematical Society, Providence, RI},
      YEAR = {2003},
     PAGES = {xvi+201},
      %ISBN = {0-8218-3440-1},
   MRCLASS = {20C30 (05E05 05E10 20C32 20P05 60C05)},
  MRNUMBER = {1984868},
MRREVIEWER = {Akihito Hora},
}

@article {kerov_okounkov_olshanski,
    AUTHOR = {Kerov, Sergei and Okounkov, Andrei and Olshanski, Grigori},
     TITLE = {The boundary of the {Y}oung graph with {J}ack edge
              multiplicities},
   JOURNAL = {Internat. Math. Res. Notices},
  FJOURNAL = {International Mathematics Research Notices},
      YEAR = {1998},
    NUMBER = {4},
     PAGES = {173--199},
      %ISSN = {1073-7928},
   MRCLASS = {05E10 (20C30 31C35 33C80)},
  MRNUMBER = {1609628},
MRREVIEWER = {Tim H. Baker},
}

@incollection {versh_ker_85,
    AUTHOR = {Kerov, S. V. and Vershik, A. M.},
     TITLE = {Locally semisimple algebras. {C}ombinatorial theory and the
              {$K_0$}-functor},
 BOOKTITLE = {Current problems in mathematics. {N}ewest results, {V}ol. 26},
    SERIES = {Itogi Nauki i Tekhniki},
     PAGES = {3--56},
 PUBLISHER = {Akad. Nauk SSSR, Vsesoyuz. Inst. Nauchn. i Tekhn. Inform.,
              Moscow},
      YEAR = {1985},
   %MRCLASS = {22D25 (18F25 19K14 46L80)},
  %MRNUMBER = {849784},
%MRREVIEWER = {\Dbar \cftil{o} Ng\d{o}c Di\cfudot{e}p},
}

@article {versh_ker_81,
    AUTHOR = {Kerov, S. V. and Vershik, A. M. },
     TITLE = {Asymptotic theory of the characters of a symmetric group},
   JOURNAL = {Funktsional. Anal. i Prilozhen.},
  FJOURNAL = {Akademiya Nauk SSSR. Funktsional\cprime ny\u{\i} Analiz i ego
              Prilozheniya},
    VOLUME = {15},
      YEAR = {1981},
    NUMBER = {4},
     PAGES = {15--27, 96},
      %ISSN = {0374-1990},
   MRCLASS = {22D10 (22D20 46L99)},
  MRNUMBER = {639197},
MRREVIEWER = {G. L. Litvinov},
}

@article {boyer,
    AUTHOR = {Boyer, Robert P.},
     TITLE = {Infinite traces of {AF}-algebras and characters of {${\rm
              U}(\infty )$}},
   JOURNAL = {J. Operator Theory},
  FJOURNAL = {Journal of Operator Theory},
    VOLUME = {9},
      YEAR = {1983},
    NUMBER = {2},
     PAGES = {205--236},
      %ISSN = {0379-4024},
   MRCLASS = {46L50 (22E45 22E65)},
  MRNUMBER = {703808},
MRREVIEWER = {G. A. Elliott},
}

@article {boyer_symplectic,
    AUTHOR = {Boyer, Robert P.},
     TITLE = {Characters of the infinite symplectic group---a {R}iesz ring
              approach},
   JOURNAL = {J. Funct. Anal.},
  FJOURNAL = {Journal of Functional Analysis},
    VOLUME = {70},
      YEAR = {1987},
    NUMBER = {2},
     PAGES = {357--387},
      %ISSN = {0022-1236},
   %MRCLASS = {22D10 (19K99 22E65 46L80)},
 % MRNUMBER = {874061},
%MRREVIEWER = {R. J. Plymen},
      % DOI = {10.1016/0022-1236(87)90117-0},
       %URL = {https://doi.org/10.1016/0022-1236(87)90117-0},
}

@ARTICLE{safonkin20,
    author = {Safonkin, N. A.},
    title = {Semifinite harmonic functions on the Gnedin--Kingman graph},
    %JOURNAL = {Zap. Nauchn. Sem. POMI},
    journal = {Zapiski Nauchnykh Seminarov POMI, Representation theory, dynamical systems, combinatorial methods. Part~XXXI},
    volume = {498},
    year = {2020},
    pages = {38--54 (Russian)},
    addendum = {English translation \fullcite{safonkin20NY}}
    %archivePrefix = {arXiv},
    %eprint = {2103.02257},
    %primaryClass = {math.CO},
}

@book {stanley1,
    AUTHOR = {Stanley, Richard P.},
     TITLE = {Enumerative combinatorics. {V}olume 1},
    SERIES = {Cambridge Studies in Advanced Mathematics},
   EDITION = {Second edition},
 PUBLISHER = {Cambridge University Press, Cambridge},
 VOLUME = {49},
      YEAR = {2012},
     PAGES = {xiv+626},
      %ISBN = {978-1-107-60262-5},
   MRCLASS = {05-02 (05A15 06-02)},
  MRNUMBER = {2868112},
}

@misc{matveev2018macdonaldpositive,
      title={Macdonald-positive specializations of the algebra of symmetric functions: Proof of the Kerov conjecture}, 
      author={Konstantin Matveev},
      year={2018},
      eprint={1711.06939},
      archivePrefix={arXiv},
      primaryClass={math.RT}
}

@book {blackadar,
    AUTHOR = {Blackadar, B.},
     TITLE = {Operator algebras},
    SERIES = {Encyclopaedia of Mathematical Sciences},
    VOLUME = {122},
      NOTE = {Theory of $C^*$-algebras and von Neumann algebras,
              Operator Algebras and Non-commutative Geometry, III},
 PUBLISHER = {Springer-Verlag, Berlin},
      YEAR = {2006},
     PAGES = {xx+517},
      %ISBN = {978-3-540-28486-4; 3-540-28486-9},
   MRCLASS = {46L05 (46L10 46L80)},
  MRNUMBER = {2188261},
MRREVIEWER = {Paul Jolissaint},
       %DOI = {10.1007/3-540-28517-2},
       %URL = {https://doi.org/10.1007/3-540-28517-2},
}

@article {versh_nikitin2010,
    AUTHOR = {Vershik, A. M. and Nikitin, P. P.},
     TITLE = {Description of the characters and factor representations of
              the infinite symmetric inverse semigroup},
   JOURNAL = {Funktsional. Anal. i Prilozhen.},
  FJOURNAL = {Funktsional\cprime ny\u{\i} Analiz i ego Prilozheniya},
    VOLUME = {45},
      YEAR = {2011},
    NUMBER = {1},
     PAGES = {16--30},
      ISSN = {0374-1990},
   %MRCLASS = {20M18 (20C32 20M30)},
  %MRNUMBER = {2848737},
%MRREVIEWER = {D. C. Jackson},
       %DOI = {10.1007/s10688-011-0002-0},
       %URL = {https://doi.org/10.1007/s10688-011-0002-0},
}
\end{document}